\documentclass[10pt]{article}

\usepackage{amsmath,makeidx,amssymb,amsthm,amscd,a4}
\usepackage{graphicx,epsfig,latexsym,index}

\newtheorem{theorem}{Theorem}[section]

\newtheorem{lemma}[theorem]{Lemma}

\newtheorem{proposition}[theorem]{Proposition}
\newtheorem{remark}[theorem]{Remark}
\newtheorem{idprobl}[theorem]{Identification Problem}

\newcommand{\D}{\protect\displaystyle}
\newcommand{\T}{\protect\textstyle}
\newcommand{\ipl}{\langle} 
\newcommand{\ipr}{\rangle} 
\newcommand{\ve}{\varepsilon}
\def\Ga{{\cal G}_\alpha}
\def\bv{{\tt{BV}}}
\usepackage{subfigure}
\newcommand{\Rr}{{\bf{\Omega}}}
\newcommand{\bydef}{\stackrel{\mbox{\scriptsize def}}{=}}
\newcommand{\QED}{\hfill \mbox{$\square$}}

\begin{document}

\title{On inverse doping profile problems for the stationary
voltage-current map}

\author{A.\,Leit\~ao${}^1$ \quad
        P.A.\,Markowich${}^2$ \quad
        J.P.\,Zubelli${}^3$ \\[2ex]}

\date{\normalsize\today}

\maketitle

\begin{center} \begin{minipage}{10cm}
${}^1$ \small Department of Mathematics, Federal University of St.\,Catarina,
\small 88040-900 Florianopolis, Brazil
({\tt aleitao@mtm.ufsc.br}) \\
${}^2$ \small Department of Mathematics, University of Vienna,
\small Boltzmanngasse 9, A-1090 Vienna, Austria
({\tt peter.markowich@univie.ac.at}) \\
${}^3$ \small IMPA, Estr.\,Dona\,Castorina 110, 22460-320 Rio de Janeiro,
Brazil ({\tt zubelli@impa.br})
\end{minipage} \end{center}
\vskip1cm

\begin{abstract}
We consider the problem of identifying possibly discontinuous doping
profiles in semiconductor devices from data obtained by\,stationary
voltage-current maps. In particular, we focus on the so-called
{\em unipolar case}, a system of PDE's derived directly from the
drift diffusion equations. The related inverse problem corresponds
to an inverse conductivity problem with partial data.
The identification issue for this inverse problem is considered. In
particular, for a discretized version of the problem, we derive a
result connected to diffusion tomography theory.
A numerical approach for the identification problem using level set
methods is presented. Our method is compared with previous results in
the literature, where Landweber-Kaczmarz type methods were used to
solve a similar problem.
\end{abstract}

%-----------------------------------------------------------------------------
\section{Introduction} \label{sec:introd}

The precise implantation of the doping profile is crucial for the desired
performance of semiconductor devices. In many applications,
there is substantial interest in replacing expensive laboratory testing by
numerical simulation and non-destructive testing, in order to minimize
manufacturing costs of semiconductors as well as for quality control.
The identification of the doping profile from indirect measurements is
called an inverse doping profile problem.
In laboratory experiments there are different types of measurement
techniques, such as {\em Laser-Beam-Induced Current} (LBIC)
\cite{FI92,FI94,FIR02}, {\em Capacitance} \cite{BELM04,BEMP01} and
{\em Current Flow} \cite{BELM04,BEM02} measurements.
These measurement techniques are related to different types of data and lead
to various inverse doping problems.
This paper is devoted to the analysis of an identification problem related
to a particular model, the so-called {\em unipolar system}, derived from
the stationary drift diffusion equations under certain simplifying
assumptions on the concentration of free carriers of positive charges
and on the recombination-generation rate. 
In this framework, the parameter function to be identified is the
{\em doping profile} . It depends on the space variables only and
represents the doping concentration, which gives the performance of
the device. It is produced by diffusion of different materials into
the silicon crystal and by implantation with an ion beam.

We shall focus on reconstruction problems based on data generated by the
{\em voltage-current} (V-C) map, i.e., an operator that takes the applied
voltage at a specified boundary part (corresponding to a semiconductor
contact) into the outflow current density on a different boundary part
(another contact). The two main contributions of this paper consist of
a theoretical identification result for a discretized version and the
analysis of a level set type method for solving the inverse doping profile
problem in the unipolar case.

The starting point of the mathematical model discussed in this paper is
the system of {\em stationary drift diffusion equations} (see system
(\ref{eq:dd-sys}) in Section~\ref{sec:model}).
This system of equations, derived more than fifty year ago \cite{vRo50},
is the most widely used to describe semiconductor devices and represents
an accurate compromise between efficient numerical solvability of the
mathematical model and realistic description of the underlying physics
\cite{Ma86,MRS90,Se84}.

This paper is organized as follows:
In Section~\ref{sec:model} we briefly introduce the {\em drift diffusion}
equations, the V-C map, and the {\em stationary linearized unipolar system}.
The latter models the direct problem related to the inverse doping profile
problem analyzed in this paper. \\
In Section~\ref{sec:idp-ident} we treat the identification issue for this
inverse problem. We do not have, at present, a theoretical result showing
uniqueness in the identification of the doping profile. However, we do
present two lines of reasoning that support the conjecture of an
identifiability result for the doping profile: The first one is based on
recent results due to Bukhgeim and Uhlmann \cite{BU02} on global uniqueness
for the local Dirichlet-to-Neumann map; The second one concerns a discretized
version of the problem that falls within the scope of tomography in the
presence of diffusion and scattering \cite{Gr92,GZ92}. \\
In Section~\ref{sec:numeric} we use a level set type method to reconstruct the
doping profile function. In this approach, a single pair of voltage-current
data is used. We compare our results with the competing Landweber-Kaczmarz
method used in \cite{BELM04} to solve a similar problem.
An analytical result concerning stability, convergence and well-posedness
of this level set method is also presented.
Section~\ref{sec:concl} is devoted to final comments and conclusions.

%-----------------------------------------------------------------------------
\section{Inverse doping profile problems} \label{sec:model}

\subsection{The semiconductor equations} \label{ssec:mod1}

The {\em drift diffusion} system of equations is the most widely used
model to describe semiconductor devices.
The mathematical modeling of semiconductor equations has developed
significantly, together with their manufacturing. The {\em basic
semiconductor device equations} were first presented, in the level
of completeness discussed here, by W.R.~van Roosbroeck in \cite{vRo50}.
Since then, it has been subject of intensive mathematical and numerical
investigation. Recent detailed expositions of the subject of modeling,
analysis and simulation of semiconductor equations can be found in
\cite{Ma86,MRS90,Se84} to cite a few.

For the sake of simplicity, we formulate the drift diffusion equations in
terms of the {\em slotboom variables} $(u,v)$.
Using an adequate change of variables, motivated by the Einstein relations,
the functions $u$ and $v$ are obtained from the electron density function
and from the hole density function respectively.
The details concerning the derivation of the model below can be found in
\cite{BELM04,BEMP01}.

The stationary drift diffusion equations consist of the Poisson equation
(\ref{eq:dd-sys1}) for the (rescaled) electrostatic potential $V$ and the
continuity equations (\ref{eq:dd-sys2}) and (\ref{eq:dd-sys3}).
\begin{subequations} \label{eq:dd-sys} \begin{eqnarray}
\lambda^2 \, \Delta V & \hskip-0.8cm \label{eq:dd-sys1}
  = \ \delta^2 \big(e^V u - e^{-V} v \big) - C & {\rm in}\ \Omega \\
{\rm div}\, J_n & \hskip-0.4cm \label{eq:dd-sys2}
  = \ \delta^4 \, Q(V,u,v,x) \, (u v - 1) & {\rm in}\ \Omega \\
{\rm div}\, J_p & \hskip-0.05cm \label{eq:dd-sys3}
  = \ - \delta^4 \, Q(V,u,v,x) \, (u v - 1) & {\rm in}\ \Omega \\
V & \hskip-1.9cm = \ V_D \ := \ U + V_{\rm bi}  & \rm on \ \partial\Omega_D
\label{eq:dd-sys4} \\
u & \hskip-2.55cm = \ u_D \ := \ e^{-U}         & \rm on \ \partial\Omega_D
\label{eq:dd-sys5} \\
v & \hskip-2.85cm = \ v_D \ := \ e^{U}          & \rm on \ \partial\Omega_D
\label{eq:dd-sys6} \\
\nabla V \cdot \nu & \hskip-1cm = \ J_n \cdot \nu \ = \
J_p \cdot \nu \ = \ 0 & \rm on \ \partial\Omega_N \label{eq:dd-sys7}
\end{eqnarray} \end{subequations}
where the densities of the electron and hole current $J_n$ and $J_p$
satisfy the current relations:
$$ J_n = \mu_n q n_i e^V \nabla u \quad {\rm and} \quad
   J_p = -\mu_p q n_i e^{-V} \nabla v \, . $$
Here the positive constants $q$ and $n_i$ denote the elementary charge
and the intrinsic charge density respectively. Moreover, $\mu_n$ and $\mu_p$
represent the (rescaled) mobilities of electrons and holes respectively.

The domain $\Omega \subset \mathbb{R}^d$ ($d=2,3$) represents the
semiconductor device. Two dimensionless positive parameters occur, namely
$\lambda$ and $\delta$, both small in many practical applications.
The function $Q$ is defined implicitly by the recombination-generation
rate function.
As far as boundary conditions are concerned, the function $U$ is the applied
potential and $V_{\rm bi}(x) := U_T \, \ln \big( n_D(x) / n_i \big)$,
where $U_T$ is the thermal voltage.

The function $C = C(x)$ denotes the doping concentration, which is
produced by diffusion of different materials into the silicon crystal and
by implantation with an ion beam.
In many technological applications, the {\em doping profile} $C$ is the
parameter that has to be identified.
Because of inaccuracies in the manufacturing process, semiconductor devices
should pass through some testing to ensure high quality.
The inverse problem we are concerned with
is related to a non destructive identification procedure, based on
experiments modeled by the {\em voltage to current} operator.

In the sequel we briefly discuss the boundary conditions
(\ref{eq:dd-sys4})--(\ref{eq:dd-sys7}).
The boundary of $\Omega$ is assumed to be divided in two nonempty parts:
$\partial\Omega = \partial\Omega_N \cup \partial\Omega_D$.
The segments of $\partial\Omega_D$ correspond to the semiconductor contacts,
where Dirichlet boundary conditions are prescribed. Differences in $U$
between different parts of $\partial\Omega_D$ correspond to the applied
bias between these two contacts.
The Neumann part of the boundary $\partial\Omega_N = \partial\Omega 
- \partial\Omega_D$ models insulating or artificial surfaces. Therefore,
a zero current flow and a zero electric field in the normal direction are
prescribed.

Existence (in weak sense) and some uniqueness results for system
(\ref{eq:dd-sys}) can be found in \cite{Ma86,MRS90}.
Under suitable regularity assumptions on the boundary conditions $u_D$,
$v_D$, $U$ and on the doping profile $C$, one can prove that system
(\ref{eq:dd-sys}) admits a weak solution $(V,u,v)$ in $(H^1(\Omega) \cap
L^\infty(\Omega))^3$. See \cite[Theorem~3.3.16]{MRS90} and
\cite[Theorem~4.2]{BEMP01}. Stronger existence results for $(H^2(\Omega)
\cap L^\infty(\Omega))^3$ can be found in \cite{Ma86}.

\subsection{The inverse problem} \label{ssec:mod2}

\begin{figure}[t]
\centerline{ \epsfxsize8cm \epsfbox{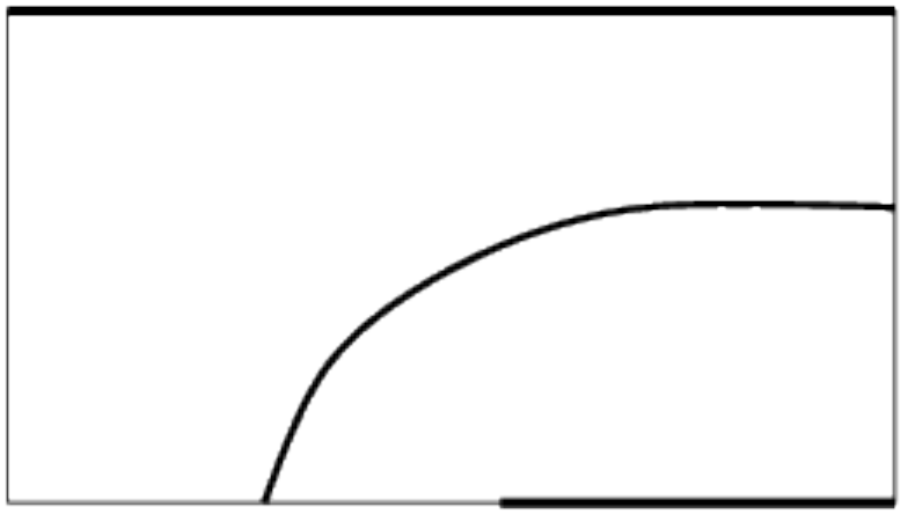} }
\caption{\small The domain $\Omega \subset \mathbb R^2$ represents a P-N
diode. The P-region corresponds to the subregion of $\Omega$, where $C < 0$.
In the N-region $C > 0$ holds. The curve between these regions is called
{\em P-N junction.} } \label{fig:diode}
\vskip-5.8cm \unitlength1cm
\centerline{
\begin{picture}(13,6)
%\put(0,0){\dashbox{0.1}(13,7){}}
\put(1.7,4.0){$\partial\Omega_N$} \put(10.6,4.0){$\partial\Omega_N$}
\put(3.4,2.2){$\partial\Omega_N$}
\put(8.0,2.2){$\Gamma_1 \subset \partial\Omega_D$}
\put(5.7,5.8){$\Gamma_0 \subset \partial\Omega_D$}
\put(3.2,4.8){\bf N-region \rm ($C > 0$)}
\put(6.7,3.4){\bf P-region \rm ($C < 0$)}
\end{picture} }
\end{figure}
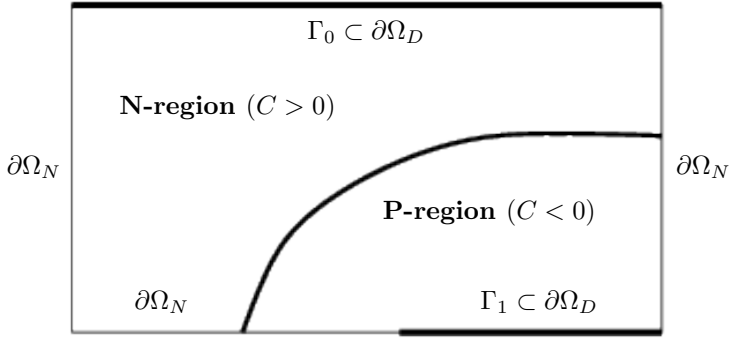

We start the discussion by introducing the {\em voltage-current} (V-C) map:
$$ \begin{array}{rcl}
   \Sigma_C: H^{3/2}(\partial\Omega_D) & \to & L^2(\Gamma_1) \\
   U & \mapsto & J\cdot\nu|_{\Gamma_1}\ =\ (J_n+J_p)\cdot\nu |_{\Gamma_1}\, ,
   \end{array} $$
where $\Gamma_1 \subset \partial\Omega_D$ corresponds to the part of the
boundary (a contact) where measurements are taken.
Notice that the map $\Sigma_C$ takes the applied voltage $U$ into the
corresponding current density.
In the inverse problem considered in this paper, the linearized V-C map
at $U = 0$ plays a key role, as we shall see later on this section.

Since the potential can be shifted by a constant, we shall assume without
lost of generalization that $U(x)|_{\Gamma_1} = 0$. In practical applications,
the applied potential $U \in H^{3/2} (\partial\Omega_D)$ is assumed to be
piecewise constant in the contacts.
To illustrate, a very simple semiconductor device, is shown in Figure~%
\ref{fig:diode}.

In the next lemma we briefly review some properties of the nonlinear operator
$\Sigma_C$. A complete proof can be found in \cite{BEMP01}.

\begin{lemma} \label{lem:vc-wd}
The current $\Sigma_C(U) \in L^2(\Gamma_1)$ is uniquely defined for
each voltage $U \in H^{3/2}(\partial\Omega_D)$ in the neighborhood of
$U = 0$, i.e., the operator $\Sigma_C$ is well-defined in the
neighborhood of $U = 0$. Moreover, $\Sigma_C$ is continuous and
continuously Fr\'echet differentiable in the neighborhood of $U = 0$.
\end{lemma}

If $U = 0$, the solution of (\ref{eq:dd-sys}) is given by $(V,u,v)
= (V^0,1,1)$, where $V^0$ is a solution of the {\em Poisson equation at
equilibrium}
\begin{equation} \label{eq:poiss-equil}
\left\{ \begin{array}{rcll}
   \lambda^2 \Delta V^0 & = & \delta^2 ( e^{V^0} - e^{-V^0} ) - C
   & {\rm in}\ \Omega \\
   V^0 & = & V_{\rm bi} & {\rm on}\ \partial\Omega_D \\
   \nabla V^0 \cdot \nu & = & 0 & {\rm on}\ \partial\Omega_N \, .
\end{array} \right.
\end{equation}
From now on, the following simplifying assumptions are made:

{\it A1)} The concentration of holes satisfy $v = 0$;

{\it A2)} No recombination-generation rate is present, i.e., $Q = 0$;

{\it A3)} The electron mobility is constant ($\mu_n = 1$) and $q=1$.

\noindent Under these assumptions, we conclude that the Gateaux derivative
of $\Sigma_C$ at $U=0$ in the direction $h \in H^{3/2} (\partial\Omega_D)$
is given by
$$ \Sigma'_C(0) h = e^{V^0} \hat{u}_\nu |_{\Gamma_1}, $$
where $\hat{u}$ and $V^0$ solve
\begin{equation} \label{eq:unipolar}
\left\{ \begin{array}{r@{\ }c@{\ }l@{\ \ }l}
   {\rm div}\, (e^{V^0} \nabla \hat{u}) & = & 0 &  {\rm in}\ \Omega \\
   \hat{u} & = & h & {\rm on}\ \partial\Omega_D \\
   J_n \cdot \nu & = & 0 & {\rm on}\ \partial\Omega_N
\end{array} \right.
\hskip0.8cm
\left\{ \begin{array}{r@{\ }c@{\ }l@{\ \ }l}
   \lambda^2 \Delta V^0 & = & e^{V^0} - C & {\rm in}\ \Omega \\
   V^0 & = & V_{\rm bi} & {\rm on}\ \partial\Omega_D \\
   \nabla V^0 \cdot \nu & = & 0 & {\rm on}\ \partial\Omega_N \\
\end{array} \right.
\end{equation}

The decoupled system (\ref{eq:unipolar}) is called {\em stationary linearized
unipolar case (close to equilibrium)}.
The inverse problem of identifying the doping profile in system
(\ref{eq:unipolar}) corresponds to the identification of $C(x)$ from
the {\em parameter-to-output map}
$$ F:
   \begin{array}[t]{rcl}
     D(F) \subset L^2(\Omega) & \to &
     {\cal L}(H^{3/2}(\partial\Omega_D);H^{1/2}(\Gamma_1)) \\
     C & \mapsto & \Sigma'_C(0)
   \end{array} $$
Since $\mu_n=1$ and $V = V_{\rm bi}(x)$ is known at $\partial\Omega_D$, the
measured current data $J_n \cdot \nu = \mu_n e^{V^0} \hat{u}_\nu$ at
$\Gamma_1$ can be directly replaced by the Neumann data $\hat{u}_\nu$.
Therefore, the inverse problem can be divided into two distinct steps:

\begin{idprobl}[Stationary linearized unipolar case] \label{probl:unipol}
\mbox{}

\begin{enumerate}
\item[(1)] Define $\gamma := e^{V^0}$ and identify $\gamma$ from the 
Dirichlet to Neumann (DtN) map
$\Lambda_\gamma: U \mapsto \gamma \hat{u}_\nu |_{\Gamma_1}$, where $\hat{u}$
solves
$$ {\rm div} (\gamma \nabla \hat{u}) \, = \, 0 \ \ {\rm in }\ \Omega\, ,\quad
   \hat{u} \, = \, U \ \ {\rm on }\ \partial\Omega_D\, ,\quad
   \hat{u}_\nu \, = \, 0 \ \ {\rm on }\ \partial\Omega_N ; $$
\item[(2)] Obtain the doping profile $C$ from: \ $C(x) = \gamma(x) -
\lambda^2 \Delta \, (\ln \gamma(x))$, $x \in \Omega$.
\end{enumerate}
\end{idprobl}

The evaluation of $C$ from $\gamma$ can be explicitely performed
(a direct problem) and is a standard procedure.
The identification issue in Problem~\ref{probl:unipol}~(1) corresponds
to the {\em electrical impedance tomography} in elliptic equations with
mixed boundary data.
For the case of the full DtN operator, i.e., $\Gamma_1 = \partial\Omega_D
= \partial\Omega$, this inverse problem has been intensively analyzed in
the literature (see, e.g., \cite{Bo02,Is98} for a survey).

%-----------------------------------------------------------------------------
\section{Inverse doping profile: Identification issue} \label{sec:idp-ident}

In this section we consider some theoretical aspects of the inverse doping
profile problem. Despite the encouraging numerical results of Section~%
\ref{sec:numeric}, at present, we do not have a theoretical result showing
uniqueness of the doping profile from V-C data measured on distinct
sub-domains of the boundary.
In Subsection~\ref{ssec:uhlmann} we present the state of the art that comes
closest to the identifiability question related to Problem~%
\ref{probl:unipol}~(1). This approach is based on recent results due to
Bukhgeim and Uhlmann~\cite{BU02} on global uniqueness for the local
Dirichlet-to-Neumann map.
In the last subsection we present a reasoning that support the conjecture
of an identifiability result for the doping profile. It concerns a discretized
version of the problem that falls within the scope of identifying the
potential of a discretized Schr\"{o}dinger equation using external
measurements. We treat this problem using techniques from the so-called
isotropic case of diffuse tomography \cite{Gr92,GZ92}.

\subsection{Global uniqueness approach} \label{ssec:uhlmann}

In the sequel we consider $\Omega$ to be 2-dimensional, unless stated
otherwise. Therefore, each current measurement is given by a function of
one space variable defined on $\Gamma_1 \subset \partial\Omega$.
Obviously, a single measurement is not sufficient to identify the doping
profile $C: \Omega\subset\mathbb R^2 \to \mathbb R$. However, adapting
some results from \cite{Na96}, related to electrical impedance tomography,
we argue in the full data case that the knowledge of the operator $F$ in
Subsection~\ref{ssec:mod2} is enough to determine $C$ uniquely.

We reason as follows: Let $V^0$ be the solution of the Poisson equation at
equilibrium in (\ref{eq:unipolar}). Given an input voltage $U \in H^{3/2}
(\partial\Omega_D)$, the output current can be identified (after rescaling)
with the Neumann data of $u$ at $\Gamma_1$, i.e., $u_\nu |_{\Gamma_1} =
\Lambda_C(U)$, where $u$ is the solution of the elliptic equation in
(\ref{eq:unipolar}). From standard results in elliptic theory, one concludes
that for a domain $\Omega$ with Lipschitz boundary, there is a one to one
relation between the solutions $V^0 \in H^2(\Omega)$ of the Poisson equation
and the potentials $C \in L^2(\Omega)$. Therefore, it is enough to consider
the problem of identifying the potential $V^0$ in (\ref{eq:unipolar}) or,
equivalently, the conductivity $\gamma = e^{V^0}$ as stated in Problem~%
\ref{probl:unipol}.

The problem of identifying conductivities from the DtN map was analyzed
by Nachman in \cite{Na96}. Adapting his result to Identification
Problem~\ref{probl:unipol} one can prove that for a bounded
$\Omega \subset \mathbb R^2$ with Lipschitz boundary, $\Gamma_1 =
\partial\Omega_D = \partial\Omega$ and $C_1$, $C_2$, two doping
profiles such that the corresponding conductivities satisfy
$$ \gamma_1, \gamma_2 \in D(F) := \{ \gamma \in W^{2,p}(\Omega), \ p > 1 ; \,
   \gamma_+ \ge \gamma(x) \ge \gamma_- > 0 \mbox{ a.e. in } \Omega \} \, , $$
the equality $\Lambda_{\gamma_1} = \Lambda_{\gamma_2}$ implies $C_1 = C_2$.

This result of Nachman has been recently improved by Astala and P\"aiv\"arinta
\cite{AP05}, who proved that any $L^\infty$ conductivity in two dimensions
can be determined uniquely from the DtN map.

We address yet another identification result (for the inverse doping profile
problem) based on the global uniqueness approach. Concerning uniqueness
results for the DtN operator with partial boundary data, this result
corresponds to the state of the art.
Let $\Omega \subset \mathbb R^n$, with $n \ge 3$, be a bounded domain with
$C^2$ boundary.
Further, let $\xi \in \mathbb R^n$ with $\|\xi\| = 1$ and $\varepsilon > 0$
be given. We define
$$ \Gamma_0 := \{ x \in \partial\Omega ; \ \ipl \nu(x), \xi \ipr >
               -\varepsilon \} , \quad
   \Gamma_1 := \{ x \in \partial\Omega ; \ \ipl \nu(x), \xi \ipr <
                \varepsilon \} $$
where $\nu(x)$ is the unit normal vector at $x \in \partial\Omega$ (notice
that $\Gamma_0 \cap \Gamma_1 \neq \emptyset$). Moreover, let $C_1$, $C_2$
be doping profiles such that the corresponding conductivities satisfy
$\gamma_1, \gamma_2 \in C^2(\overline{\Omega})$ and $\gamma_j(x) \ge \gamma_-
> 0$ a.e. in $\Omega$, $j=1,2$.
Then, the equality $\Lambda_{\gamma_1} = \Lambda_{\gamma_2}$ implies
$C_1 = C_2$ (see \cite{BU02}).

Notice that this result applies to 3-dimensional domains $\Omega$ with
regular boundary and, moreover, $\partial\Omega = \partial\Omega_D =
\Gamma_0 \cup \Gamma_1$, $\Gamma_0 \cap \Gamma_1 \neq \emptyset$, i.e.,
the contacts where the voltage is prescribed ($\Gamma_0$) and where the
current is measured ($\Gamma_1$) overlap.

\subsection{The discrete Schr\"{o}dinger equation with partial DtN data}

In this section, we consider the characterization problem for the
Schr\"{o}dinger operator potential $V$ given partial information on
the Dirichlet-to-Neumann map $\Lambda^{V}$ associated to the problem
\begin{equation} \label{SSS}
\left\{ \begin{array}{l} 
           -\Delta w + V w = 0 \mbox{ in } \Omega \\
           w \big|_{\partial \Omega} = \phi
        \end{array}
\right.
\end{equation}
It is well-known that the change of variables
\begin{equation} \label{changevars}
w = \gamma^{1/2} u \quad \mbox{ and } \quad
V = \gamma^{-1/2} \Delta \gamma^{1/2}
\end{equation}
establishes a $1-1$ correspondence between the solutions of (\ref{SSS})
and those of
\begin{equation} \label{EEE}
\left\{ \begin{array}{l} 
  {\mathrm div} (\gamma \nabla u) = 0 \mbox{ in } \Omega \\
  u \big|_{\partial \Omega} = \gamma^{-1/2} \big|_{\partial \Omega} \phi
\end{array}  
\right.
\end{equation}
The Dirichlet-to-Neumann map for (\ref{SSS}) is related to that
of (\ref{EEE}) by
\begin{equation} \label{DTN}
 \Lambda^{V}(\phi) = \gamma^{-1/2} \Lambda_{\gamma} ( \gamma^{-1/2} \phi)
 + \frac{1}{2\gamma}
 \frac{\partial \gamma}{\partial n} \phi \, .
\end{equation}
It is clear that the knowledge of the DtN map $\Lambda^{V}$ for
Equation~(\ref{SSS}) is equivalent to knowledge of its counterpart
$\Lambda_{\gamma}$ for (\ref{EEE}). Furthermore, any restriction of
$\Lambda^{V}$ to $\phi$ supported on a subset $\Gamma_{0}$ of the
boundary corresponds to the restriction of $\Lambda_{\gamma}$ supported
on this set $\Gamma_{0}$.
If we consider current measurements taken in a subset $\Gamma_{1}$
contained in $\partial \Omega$, then, at the level of $\Lambda^{V}$
this means that we will only consider the information from $\Lambda^{V}$
on $\Gamma_{1}$. Let us call such map $\Lambda^{V}\big|_{\Gamma_{0},
\Gamma_{1}}$.

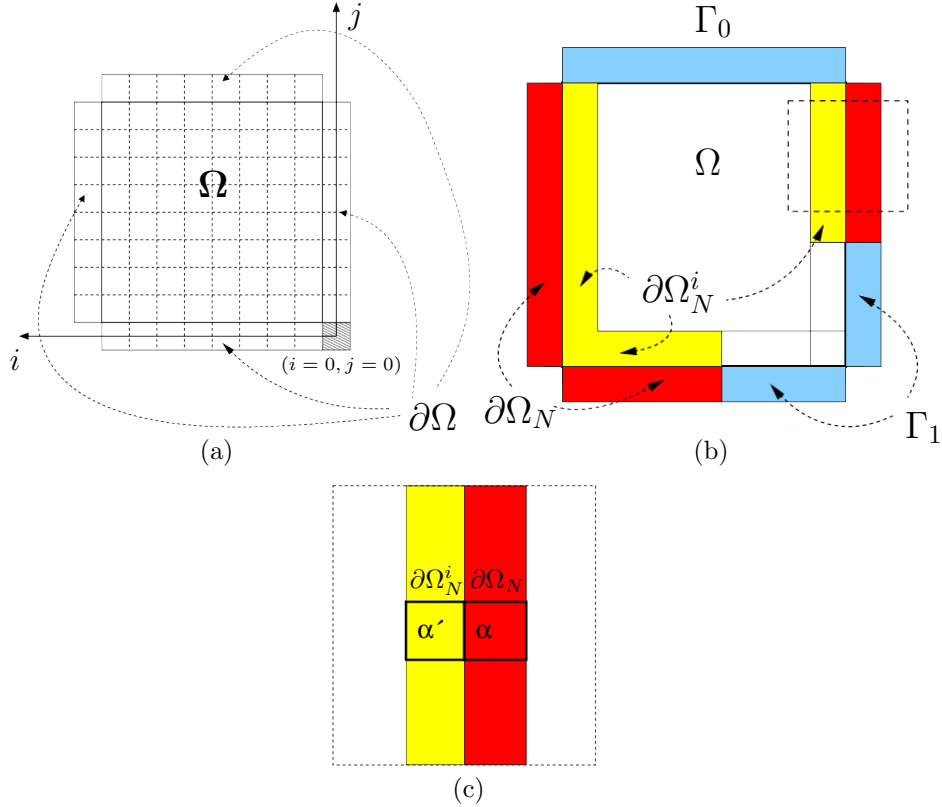
\begin{figure}[ht] %JPZ figure
\centerline{ \epsfxsize6.0cm \epsfbox{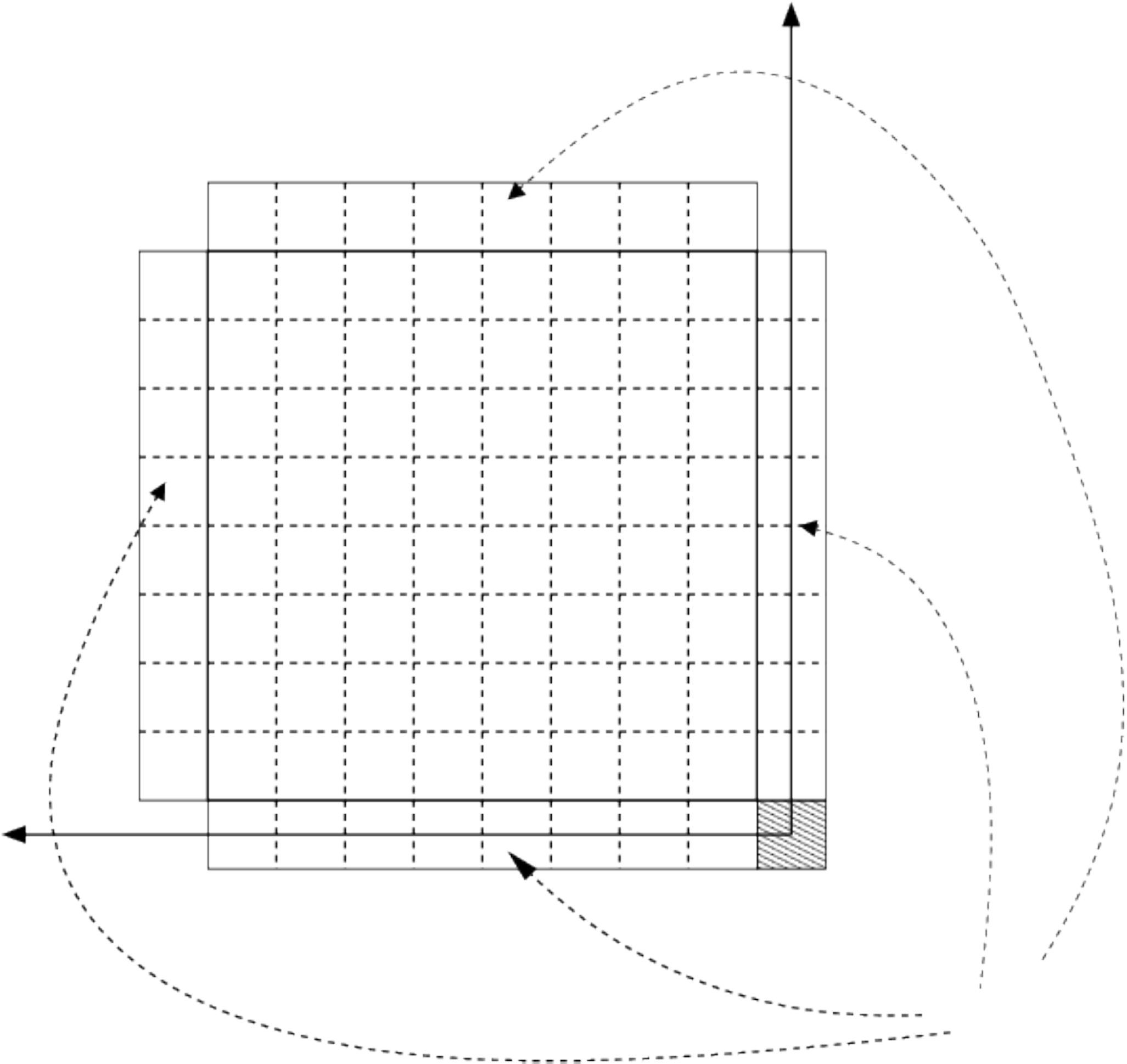} \hfil
             \epsfxsize5.6cm \epsfbox{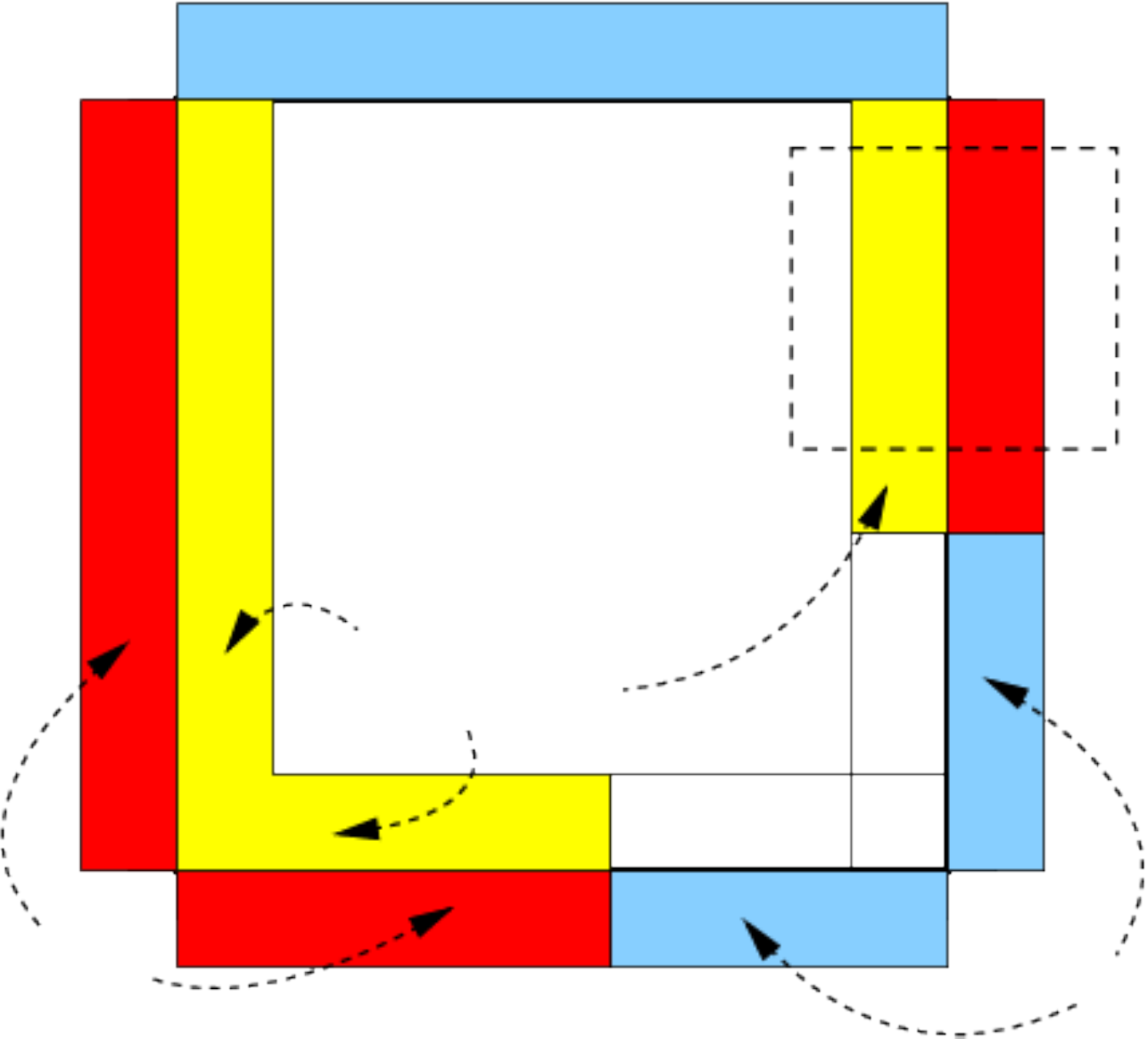} }
\centerline{\hfil (a) \hskip6.0cm  (b) \hfil}
\medskip
\centerline{\epsfxsize3.5cm \epsfbox{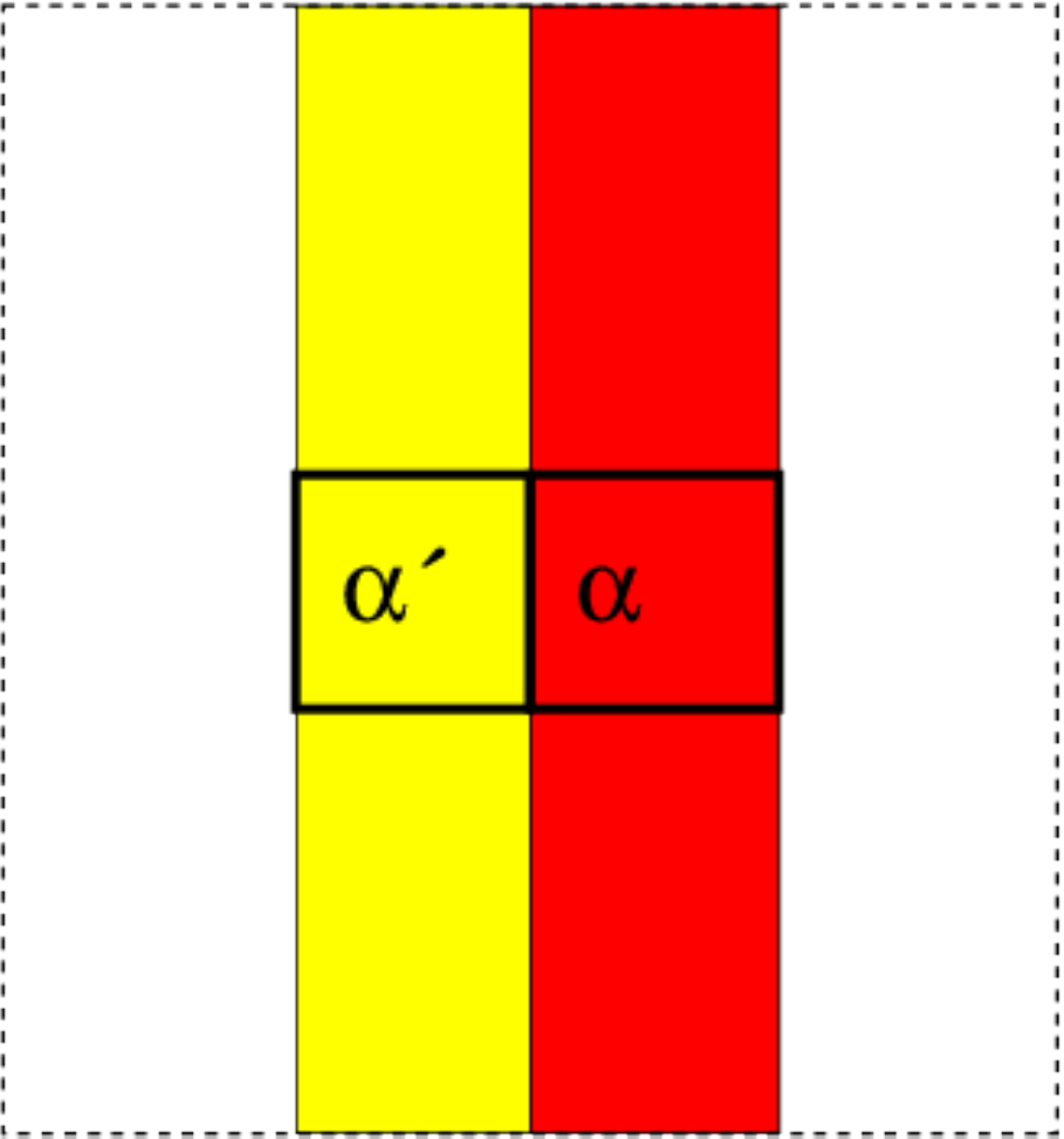} }
\centerline{(c)}
\caption{\small Picture (a) shows the discretized region $\Omega$ under
consideration and its boundary $\partial \Omega$. It also shows the origin
$(0,0)$ on the bottom right corner.
On Picture (b) the boundary parts $\partial\Omega_N$, $\partial\Omega_D =
\Gamma_0 \cup \Gamma_1$ are shown.
Picture (c) zooms in the squared region marked in (b). It shows the boundary
$\partial\Omega_N$ and it's adjacent interior part $\partial\Omega_N^i$.}
\label{discrete1}
\vskip-7.4cm \unitlength1cm
\centerline{
\begin{picture}(13,7.6)
%\put(0,0){\dashbox{0.1}(13,12.9){}}
\put(2.9,10.2){\Large $\mathbf\Omega$}
\put(5.7,7.1){\Large $\partial\Omega$}
\put(4.0,7.9){\scriptsize $(i=0,j=0)$}
\put(0.4,7.9){\large $i$} \put(4.9,12.5){\large $j$}
\put(9.5,10.5){\Large $\Omega$}
\put(9.5,12.4){\Large $\Gamma_0$} \put(12.3,7.0){\Large $\Gamma_1$}
\put(6.7,7.2){\Large $\partial\Omega_N$}
\put(8.8,8.8){\Large $\partial\Omega_N^i$}
\put(5.6,5.0){       $\partial\Omega_N^i$}
\put(6.4,5.0){       $\partial\Omega_N$}
\end{picture} }
\end{figure}

To the best of our knowledge, there is no characterization result of $V$
based on $\Lambda^{V}\big|_{(\Gamma_{0},\Gamma_{1})}$ when $\Gamma_{0} \cap
\Gamma_{1} = \emptyset$. 
We explore here the discrete analogue of the Dirichlet-to-Neumann
characterization problem with partial data for the Schr\"odinger operator.
In this context we consider a discretization $V_{ij}=V(x_i,y_j)$ of
$V:\Omega \rightarrow \mathbb{R}$ for $(i,j)\in \Rr\bydef \left\{ (i,j) |
1 \le i,j \le N, i,j \in {\mathbb{Z}} \right\}$. For a mesh size
$\Delta x = \Delta y = \epsilon$, the first equation in (\ref{SSS}) is
replaced by
\begin{equation} \label{disc1}
u_{ij} = \frac{1}{4+ \epsilon^2 V_{ij}}
         ( u_{i+1,j } + u_{i-1,j} + u_{i, j+1} + u_{i, j-1} ) 
\qquad \mbox{\rm for } (i,j) \in \Rr \, ,
\end{equation}
We define $w_{ij} = 4/(4+\epsilon^{2} V_{ij})$ and consider the set of
equations described by
\begin{equation} \label{system1}
u_{ij} - \frac{w_{ij}}{4} \left( u_{i+1,j } + u_{i-1,j} + u_{i, j+1}
       + u_{i, j-1}  \right) = 0 \, \mbox{, where $(i,j)\in \Rr$.}
\end{equation}
We remark that except for minor modifications, in what follows, we could
use $1 \le i \le N_1$ and $1 \le j \le N_2$ (see Figure~\ref{discrete1}).

The system of equations defined by (\ref{system1}) must be supplemented
with suitable boundary conditions. In \cite{Gr92,GZ92}, Dirichlet type
boundary conditions were imposed for $u_{i,j}$ whenever
$(i,j)\in \partial \Rr$, where $ \partial \Rr$ is the set of points
$(i,j)$ with $0\le i,j \le N+1$ where either $i\in \{0, N+1\}$ or
$j \in \{ 0, N+1 \}$, but not both. See Figure~\ref{discrete1}.
More precisely, in \cite{Gr92,GZ92} one imposes  the condition
\begin{equation} \label{boundary1}
u_{d} = \delta_{d} \quad \forall d=(i_0,j_0) \in \partial \Rr \, ,
\end{equation}
where $\delta_{d}(l) \bydef 1 $ if $d = l$ and $0$ otherwise.

If $0\le w_{ij} \le 1$ for all $(i,j) \in \Rr$, then the problem
(\ref{system1}) with boundary conditions (\ref{boundary1}) has a
natural probabilistic interpretation. Namely, $u_{ij}$ represents
the probability that a particle undergoing a random walk with
absorption will reach the site $d=(i_0,j_0)$ given that at each site
$\alpha = (\alpha_1,  \alpha_2)$ it has a survival probability
$w_{\alpha}$ for $\alpha \in \Rr$. See \cite{Gr92}. We remark that
a sufficient condition for $w_{ij} \in (0,1)$ is that $V_{ij}>0$.
In what follows we will rely heavily on such interpretation and the
notation presented in \cite{Gr92,GZ92}.
We shall extend some of the results therein to allow more general boundary
conditions in the identification of the doping profile.
We refer the reader to Figure~\ref{discrete1} where the different
boundary conditions are depicted. In detail, 
the boundary region $\partial \Rr$ will be decomposed into two parts,
$\partial \Rr_{N}$ and $\partial \Rr_{D}$. 
Such regions have corresponding internal adjacent regions 
$\partial \Rr_{N}^{i}$ and $\partial \Rr_{D}^{i}$. 
On $\partial \Rr_{N}$
homogeneous Neumann boundary conditions will be imposed. In this
discretized setting, this means that the values of $u_{\alpha}$ on
pixels  $\alpha \in \partial \Rr_{N}$ and on the adjacent one
$\alpha' \in \partial \Rr_{N}^{i}$ coincide. See Figure~\ref{discrete1}c. 
The region $\partial \Rr_{D}$ will be further subdivided
into two regions $\Gamma_0$ and $\Gamma_1$. On $\Gamma_0$ we will
impose nonhomogeneous Dirichlet data whereas on $\Gamma_1$ we impose
homogeneous Dirichlet data. 
Here again, for $l=0$ or $1$ we denote by  
$\Gamma_{l}^{i}$ the interior region adjacent to $\Gamma_{l}$.
The measurements correspond to normal
derivatives on $\Gamma_{1}$. In other words,
$u_{\alpha} - u_{\alpha'}$ for $\alpha \in \Gamma_{1}$ and
$\alpha' \in \Gamma_{1}^{i}$ with $\alpha'$ adjacent to $\alpha$.
Since $u_{\alpha} = 0$ for $\alpha \in \Gamma_{1}$ this corresponds
to evaluating $u_{\alpha'}$ for $\alpha' \in \Gamma_{1}^{i}$.

The first natural question to be addressed is the well posedness of
the direct problem. It is answered by the following:

\begin{proposition} Given a distribution of values
$w=(w_{i,j})_{1\le i, j \le N} \in (0,1)^{N\times N}$ the system of
equations in (\ref{system1}) endowed with the boundary conditions
\begin{eqnarray}
u_{\alpha} & = & u_{\alpha'}  \mbox{ for } \alpha \in \partial \Rr_{N}
\mbox{ adjacent to }
\alpha' \in \partial \Rr_{N}^{i} \label{bd1} \\
u_{\beta} & = & \delta_{d} \mbox{ for } \beta \in \Gamma_{0}\label{bd2} \\
u_{\gamma} & = & 0 \mbox{ elsewhere on } \partial \Rr \, , \label{bd3}
\end{eqnarray}
has a unique solution for each $d\in \Gamma_{0}$. Furthermore, this
solution depends rationally on the components of the array $w$.
\end{proposition}
{\it Proof:} \rm
Let us notice that we have a (sparse) system of $N^2$ equations in the
$N^2$ unknowns $((u_{ij}))$. The equations for the sites $(i,j)$ with
$2 \le i, j \le N-1$ are precisely those given by (\ref{system1}),
whereas for the sites $\alpha' = (i,j) \in \partial \Rr^{i}_{N}$ or
$\partial \Rr^{i}_{D}$ require us to use the boundary conditions.
The variables $u_{\alpha}$ in the site $\alpha$ adjacent to
$\alpha' \in \partial \Rr^{i}_{N}$ coincides with $u_{\alpha'}$.
Thus, the corresponding equation has to be modified accordingly.
On the other hand, if $\alpha' \in \partial\Rr^{i}_{D}$ then the
value of $u_{\alpha}$ must be $\delta_{d}(\alpha)$.
In the sites adjacent to the Dirichlet boundary, or in the interior
sites, the diagonal element of the matrix representing the
system~(\ref{system1}) is 1. On the sites adjacents to the Neumann
boundary the value of $w_{ij}$ must be changed to $w_{ij}/(1-(w_{ij}/4))$.
In either case, after incorporating the boundary conditions (of mixed
Neumann and Dirichlet type)  the matrix representing the problem is
strictly diagonally dominant.
Thus the sytsem of equations is uniquely solvable, and the solution
depends rationally on the coefficients $w_{ij}$.
\QED

\begin{remark}
The assumption $w_{ij} \le 1$ for all $i$ and $j$ is crucial for the
above argument. This is ensured, for example, if $V_{ij}>0$ for
all $i$ and $j$, which in turn can be guaranteed if $V(x)$ is positive.
\end{remark}

\begin{remark} \label{neumann}
The vanishing Neumann boundary conditions can be recast so as to preserve
the probabilistic interpretation of the problem as follows: Suppose that
$(i,j) \in \partial \Rr^{i}_{N}$ is adjacent to $(i-1,j) \in \partial\Rr_{N}$
(similar considerations hold at the other points $(i,j) \in \partial
\Rr^{i}_{N}$). Then, the Equation~(\ref{system1}) for this site becomes
\begin{equation}\label{system2}
u_{ij} - \frac{w_{ij}}{4- w_{ij}}
\left( u_{i+1,j } + u_{i-1,j} + u_{i, j+1} + u_{i, j-1}  \right) = 0 \, .
\end{equation}
\end{remark}

\begin{remark}
Since the variable $u_{i-1,j}$ and the coefficient $w_{i,j}$ do not
appear in any other equation in the system, we could reinterpret
Equation~(\ref{system2}) as
$
u_{ij} = (w_{ij}^{\rm eff}/3)
\left( u_{i+1,j } + u_{i-1,j} + u_{i, j+1} + u_{i, j-1}  \right)
$, 
with $w^{\rm eff}_{ij} = 3 w_{ij} /(4 - w_{ij})$. Notice that
$w^{\rm eff}_{ij} \in (0,1)$ if $w_{ij}\in (0,1)$. Thus, for all
practical purposes, the equations associated to the Neumann boundary
sites could be replaced by equivalent equations with vanishing
Dirichlet boundary conditions.
\end{remark}

\begin{figure}[t] %JPZ figure
\centerline{ \epsfxsize5.0cm \epsfbox{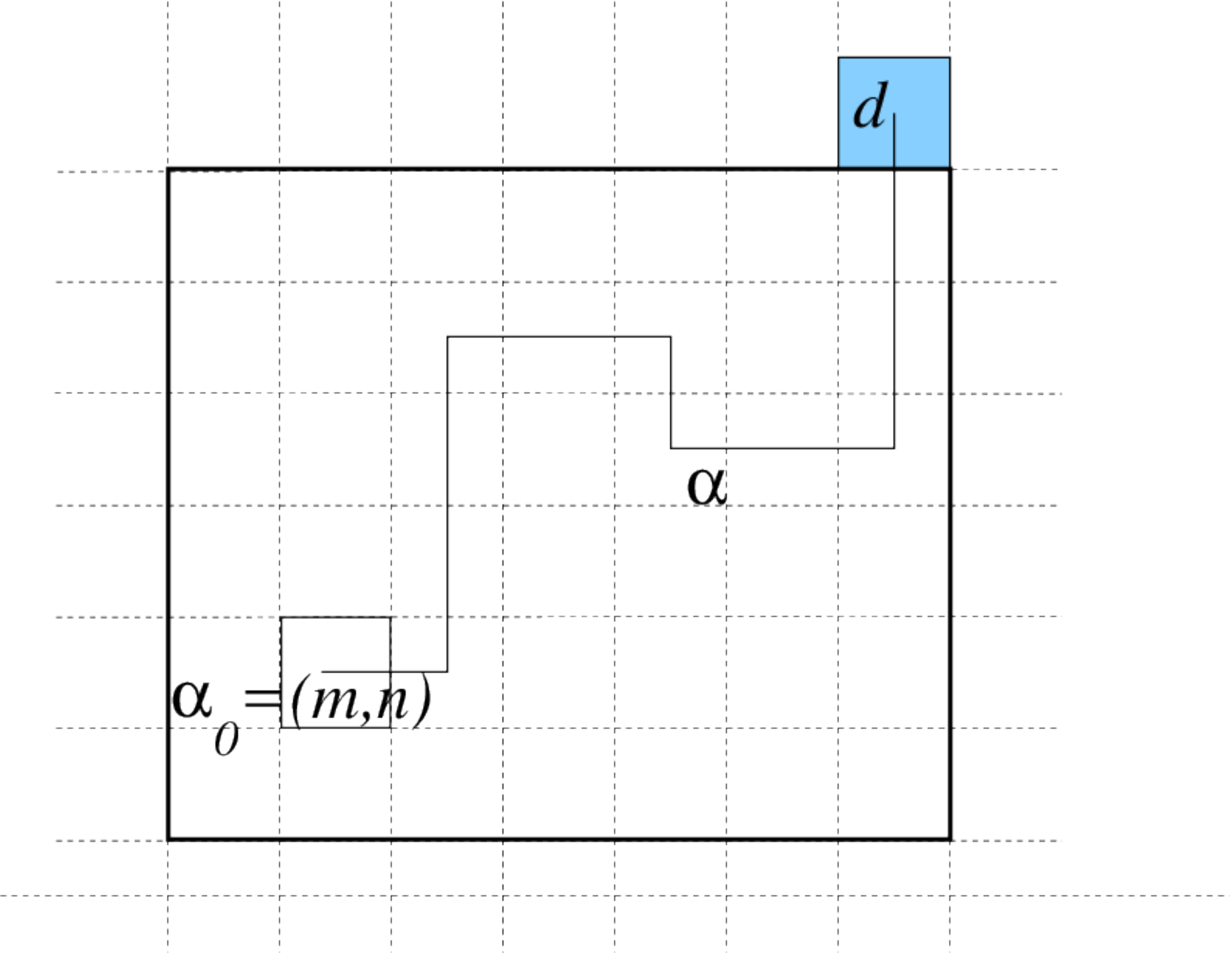} \hfil
             \epsfxsize5.0cm \epsfbox{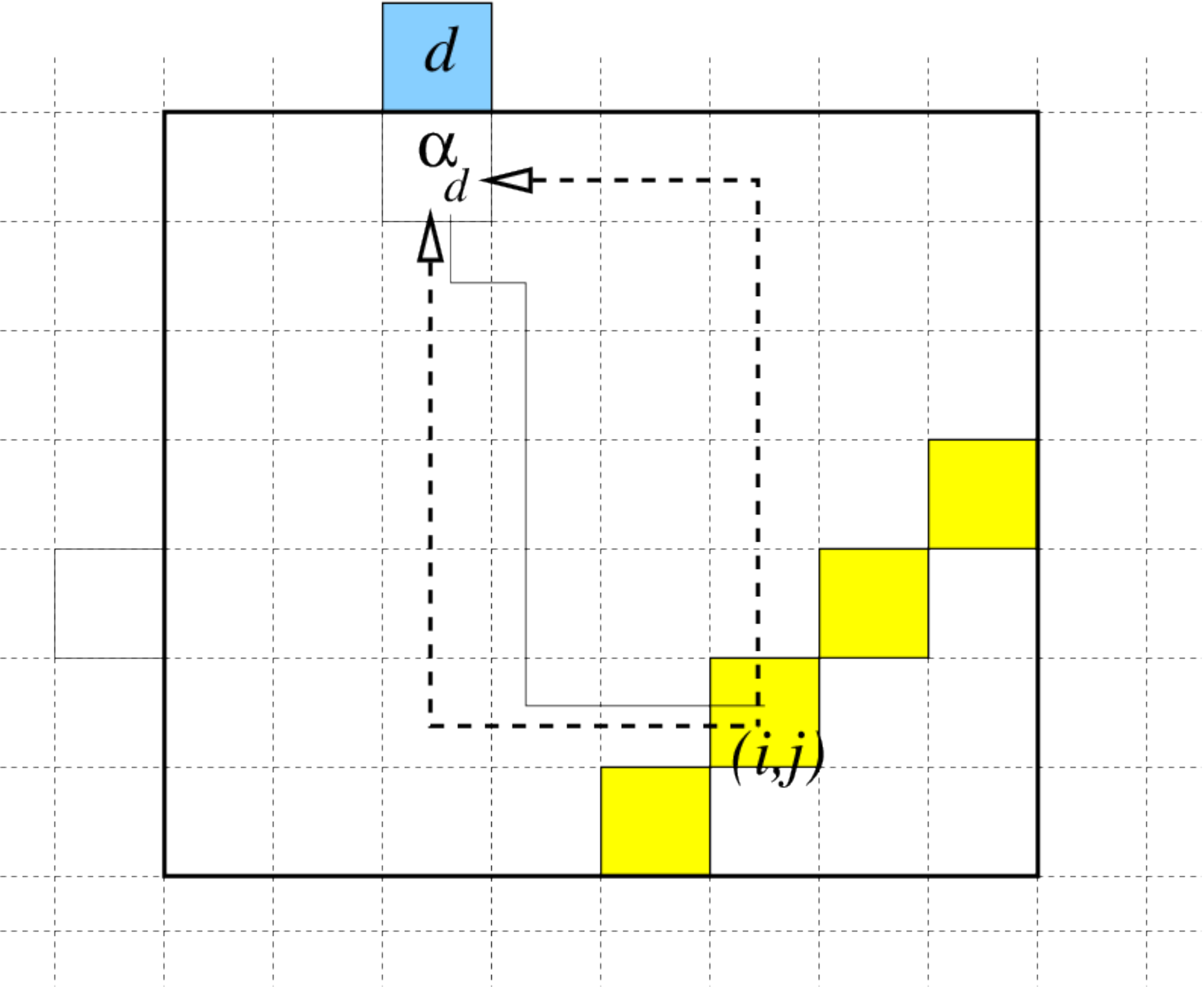} }
\centerline{\hfil (a) \hskip5.6cm  (b) \hfil}
\caption{\small Picture (a) shows an example of a path $\alpha$ connecting an
internal point to a boundary point $d$. On picture (b) an example of a few
minimal length paths connecting an internal point to the point $\alpha_d$
adjacent to a detector $d$ is shown.}
\label{discrete2}
\end{figure}

\begin{remark} \label{kac}
In \cite{Gr92} a crucial role is played by the probabilistic interpretation
of the system of equations (\ref{system1}) in the solution of the inverse
problem of the so-called isotropic diffuse tomography problem. See also
\cite{SZ2001}. In particular, the following Feynman-Kac type formula holds
for a fixed $d\in \Gamma_{0}$
$$
u_{mn} = \sum_{\alpha\in {\mathcal{P}}_{(m,n)}^{d}} \prod_{s\in \alpha}
t_{\alpha_s}^{\alpha_{s+1}} \, ,
$$
where ${\mathcal{P}}_{(m,n)}^{d}$ denotes the set of all paths connecting
the site $(m,n)$ to the boundary site $d$, and a path $\alpha$ consists of
an ordered set of successively adjacent sites starting at a neighbor to
$(m,n)$ and ending at $d$, and $t_{\alpha_s}^{\alpha_{s+1}}$ denotes the
transition probability from the site $\alpha_s$ to the site $\alpha_{s+1}$.
Thus, $t_{\alpha_0}^{\alpha_1}= w_{mn}/4 , \cdots$.
\end{remark}

We now turn our attention to the inverse problem. We define the
restricted (discrete) DtN map $\Lambda_{\Gamma_0,\Gamma_1}^{w}$ which
assigns Dirichlet data supported on $\Gamma_0$ to Neumann measurements
on $\Gamma_1$.
Our next goal is to prove an identification result that is similar in spirit
to the main result of \cite{Gr92}.
It implies that in the discrete context, and under suitable hypothesis
on the potential, one can determine such potential in the interior of a
region defined by the current measurement boundary using voltage to current
measurements.
The larger the boundary $\Gamma_1$ in Figure~2.b, the larger the region where
the potential can be uniquely determined from voltage to current measurements
(provided the total length of $\Gamma_1$ does not exceed the length of the
side of the device). More precisely, we have the following:

\begin{theorem}
For a dense open set of values $w \in (0,1)^{N\times N}$, the map
$\Lambda_{\Gamma_0,\Gamma_1}^{w}$ uniquely determines the values of
$w_{ij}$ for $(i,j)\in \Rr$ satisfying $2 \le i+j \le p'+1$, $2p' \le N+1$,
provided the support of the Dirichlet data contains the points
$(N+1,N),\cdots,(N+1,N-p'+1)$ and $p'$ is smaller than the size of
one of the sides of $\Gamma_1$.
\end{theorem}
{\it{ Proof:}}
The argument follows closely that of \cite{Gr92} by proceeding along the
diagonals. The $p$-th diagonal is defined by the sites $(i,j) \in \Rr$
such that
$ i+j = p+1 \, . $
For instance, the very first diagonal, associated to $p=1$, leads to the
equation
\begin{equation} \label{eq4}
\widehat{V}_{11} z_{11}^{d} - (z_{12}^{d} + z_{10}^{d} + z_{01}^{d}
+ z_{21}^{d} ) = 0
\end{equation}
where $(z_{ij}^{d})$ denotes the solution of system~(\ref{system1}) 
% $$ u_{ij} - \frac{w_{ij}}{4} \left( u_{i+1,j } + u_{i-1,j} + u_{i, j+1}
 %       + u_{i, j-1}  \right) = 0 \, \mbox{, where $(i,j)\in \Rr$.} $$
with boundary conditions~(\ref{bd1})--(\ref{bd3}), and Dirichlet data specified
as $\delta_{d} $ for $d\in \Gamma_0$. Furthermore, we shall use the notation
$\widehat{V}_{ij} \bydef 4/w_{ij}$.
In this simple case, we see that in Equation~(\ref{eq4})
$z_{10}^{d} = z_{01}^{d} = 0$ and $z_{11}^{d}, z_{12}^{d},z_{21}^{d}$
are all boundary measurements, and thus we can recover $\widehat{V}_{11}$.
The next diagonal ($p=2$) yields for each detector $d$:
\begin{eqnarray*} % \label{eq??}
\widehat{V}_{12} z_{12}^{d} - (z_{13}^{d} + z_{11}^{d} + z_{02}^{d} +
 z_{22}^{d} ) & = & 0 \\
\widehat{V}_{21} z_{21}^{d} - (z_{22}^{d} + z_{20}^{d} + z_{11}^{d} +
 z_{21}^{d} ) & = & 0
\end{eqnarray*}
Here the unknowns are $\widehat{V}_{12}$, $\widehat{V}_{21}$, and
$z_{22}^d$. The remaining variables, $ z_{12}^{d}$, $z_{21}^{d}$,
$z_{13}^{d}$, $z_{11}^{d}$, $z_{20}^{d}$, $z_{02}^{d}$ are all
boundary values or known from the measurements.
Upon choosing two distinct detectors we are led to the system:
\begin{equation} \label{eqstar1}
\left\{
\begin{array}{lll}
\widehat{V}_{12} z_{12}^{d_1} -
      (z_{13}^{d_1} + z_{11}^{d_1} + z_{02}^{d_1} + z_{22}^{d_1} ) & = & 0 \\
\widehat{V}_{21} z_{21}^{d_1} -
      (z_{22}^{d_1} + z_{20}^{d_1} + z_{11}^{d_1} + z_{21}^{d_1} ) & = & 0 \\
\widehat{V}_{12} z_{12}^{d_2} -
      (z_{13}^{d_2} + z_{11}^{d_2} + z_{02}^{d_2} + z_{22}^{d_2} ) & = & 0 \\
\widehat{V}_{21} z_{21}^{d_2} -
      (z_{22}^{d_2} + z_{20}^{d_2} + z_{11}^{d_2} + z_{21}^{d_2} ) & = & 0 \, ,
\end{array}
\right.
\end{equation}
where the unknown is $(\widehat{V}_{12},\widehat{V}_{21},z_{22}^{d_1},
z_{22}^{d_2})$.
The system has a unique solution iff its determinant, which is given by
$z_{12}^{d_1} z_{21}^{d_2} - z_{21}^{d_1} z_{12}^{d_2} $, does not vanish.
In this case, as a byproduct of the solution we also determine $z_{22}^{d_1}$
and $z_{22}^{d_2}$. The latter will be used in the next step, together with
a possible collection of other values of $z_{22}^{d}$ for $d\in \left\{ d_1,
d_2,\cdots,d_{\rm max} \right\}$.
In general, for $p\ge 1$, the equations associated to $i+j=p+1$ and detector
$d$ take the form
\begin{equation} \label{eqstar3}
\left\{
\begin{array}{cll}
\widehat{V}_{1,p} z_{1,p}^{d} -
      (z_{1,p+1}^{d} + z_{1,p-1}^{d} + z_{0,p}^{d} + z_{2,p}^{d} ) & = & 0 \\
\widehat{V}_{2,p-1} z_{2,p-1}^{d} -
      (z_{2,p}^{d} + z_{2,p-2}^{d} + z_{1,p-1}^{d} + z_{3,p-1}^{d} ) & = & 0 \\
\vdots & =  & \vdots \\
\widehat{V}_{p-1,2} z_{p-1,2}^{d} -
      (z_{p-1,3}^{d} + z_{p-1,1}^{d} + z_{p-2,2}^{d} + z_{p,2}^{d} ) & = & 0 \\
\widehat{V}_{p,1} z_{p,1}^{d} -
      (z_{p,2}^{d} + z_{p,0}^{d} + z_{p-1,1}^{d} + z_{p+1,1}^{d} ) & = & 0 \, .
\end{array}
\right.
\end{equation}
Notice that if we assume that the values of $z_{i',j'}^{d}$ have all been
determined (or measured) for $i'+j' \le p+1 $ then the unknowns become
$$\widehat{V}_{1,p} , \widehat{V}_{2,p-1}, \cdots, \widehat{V}_{p,1}
\mbox{  and  } z_{1,p+1}^{d}, z_{2, p}^{d}, \cdots, z_{p+1,1}^{d} \, . $$
We now order the detectors $d_1,d_2,\cdots,d_m$,  successively from left
to right, on the region $\Gamma_0$ of  Figure~\ref{discrete1}(b). By
detectors we mean positions where the Dirichlet data is taken to be
$\delta_{d}(i,j)$. \footnote{We recall that in the region $\Gamma_0$ we
control the voltages and by placing such detectors in this region we are
defining a basis for the space of input voltages.} The given data consists
of the currents in the region adjacent to $\Gamma_1$. Since on $\Gamma_1$
we have $u=0$, knowledge of the currents tantamounts to knowledge of the
values of $u_{\alpha}^{d}$ for $\alpha\in \{ 0 \} \times \{1,\cdots,p'\}$
or $\alpha \in \{1,\cdots,p' \} \times \{ 0 \}$.

We now introduce the following inductive hypothesis: \\
\noindent {\sc H1:} For a generic (open and dense) set $\mathcal{A}$
of the space of unknowns $((\widehat{V}_{ij}))\in (1,\infty)^{N\times N}$
one can solve the system of equations (\ref{eqstar3}) for the variables
$z_{ij}^{d}$ with $i+j \le p+2$, $d\in \{d_1, d_2,\cdots, d_p \}$, and
$\widehat{V}_{ij}$ for $i+j\le p+1$ in terms of the given data.

In the present context, by data we mean the values of $z_{ij}^{d}$ for
which any of the indices $i$ or $j$ belongs to the set $\{ 0,1 \}$ and
$d\in \{d_1, d_2,\cdots, d_p \}$.
The validity of the induction hypothesis for $p=1$ derives from the remark
following Equation~(\ref{eq4}) above. The inductive step relies on the fact
that in order to go from $p$ to $p+1$ we have to solve a system of equations
based on (\ref{eqstar3}) for detectors $d_1,\cdots,d_p $. This in turn, is
equivalent to showing that the determinant
\begin{equation} \label{bigdet}
D_{p} % (d_1,\cdots,d_p)
\bydef \left| \begin{array}{llllll}
z_{ 1,p}^{ d_1} & z_{1,p }^{ d_2} &
z_{1,p }^{ d_3} & \hdots & z_{ 1,p}^{ d_{p-1}} &z_{ 1,p}^{ d_p} \\
z_{ 2,p-1}^{ d_1 } & z_{  2,p-1}^{ d_2} &
z_{  2,p-1}^{ d_3} & \hdots & z_{  2,p-1}^{ d_{p-1}} &z_{ }^{ d_p} \\
\vdots & \vdots & \vdots & \vdots & \vdots & \vdots \\
z_{p-1,2 }^{ d_1} & z_{ p-1,2 }^{ d_2} &
z_{ p-1,2 }^{ d_3} & \hdots & z_{ p-1,2 }^{ d_{p-1}} &z_{ p-1,2 }^{ d_p} \\
z_{p,1 }^{ d_1} & z_{ p,1 }^{ d_2} &
z_{ p,1 }^{ d_3} & \hdots & z_{ p,1 }^{ d_{p-1}} &z_{ p,1 }^{ d_p} \\
\end{array} \right|
\end{equation}
does not vanish in the set $\mathcal{A}$. Although the technique we
employ here is the very same used in \cite{Gr92}, the crucial difference
is that in our case the determinant $D_{p}$ consists of detectors on the
opposite side from where the measurements are being taken. More precisely,
we show that the analytic function $D_{p}(w)$ is not identically zero in
a neighborhood of $w=0$.
Another difference from the situation in \cite{Gr92} is the fact that we
have Neumann type boundary conditions in part of the boundary. This, however,
causes no further difficulty at the light o% $$ \partial \Rr \bydef \left\{ (i,j) | i \in \left\{0, N+1 \right\} 
% \mbox{ or}  j  \in \left\{0, N+1 \right\} \right\} \setminus
% \left\{ (0,0), (0, N+1), (N+1,0),(N+1, N+1) \right\} \, . $$
%
f Remark~\ref{neumann}.

To complete the proof, it thus remain to show that if we take all values
of $w_{ij}=\rho$ and let $\rho \rightarrow 0$, then under the assumption
that $p\le p'$, $D_{p} = A(p) \rho^{L(p)} + \mathcal{O}(\rho^{L(p)}+1)$
with $A(p)\ne 0$ and $L(p)$ depending only on geometric parameters
associated to the size of the grid and the location of the detectors and
the diagonal $p$. To prove this claim, we start by noticing that because of
Remark~\ref{kac}, when $\rho \rightarrow 0$, we have
$z_{ij}^{d}(\rho) = A_{i,j}^{d,p} \rho^{\ell(p,i,j)+1} +
\mathcal{O}(\rho^{\ell(p,i,j)+2})$, where $\ell(p,i,j)$ is the length of
the smallest path connecting the site $(i,j)$ to the point $\alpha_d$
in $\Gamma_0$ adjacent to the detector $d$. See Figure~\ref{discrete2}(b).
Furthermore, $A_{i,j}^{d,p}$ denotes the number of paths in the region
$\Rr$ of minimal length $\ell(p,i,j)$ connecting $(i,j)$ to $\alpha_{d}$.
If we assume that the coordinates of $\alpha_{d}=(i',j')$ then it is easy
to check that $\ell(p,i,j)=|i'-i|+|j'-j|$ and that the number of such paths
is given by
\begin{equation} \label{numberpaths}
A_{i,j}^{d,p} = \left( \begin{array}{c}
                  |i'-i|+|j'-j| \\ |i'-i|
                \end{array} \right) =
\left( \begin{array}{c}
  \ell(p,i,j) \\ |i'-i|
\end{array} \right) =
\left( \begin{array}{c}
 \ell(p,i,j) \\ |j'-j|
\end{array} \right)
\end{equation}
A straightforward combinatorial argument gives that $A(p)\ne 0$ provided
$2p'\le N+1$.  \qed

The results presented in this subsection, although following the main
ideas in \cite{Gr92} lead to a much more difficult problem than that presented
therein. In particular, it is not clear how to go beyond $ p^{\prime}$. 
In fact, the hypothesis that $2p'\le N+1$ is crucial in the
above argument, and
although it seems it could be relaxed we do not have a proof of this fact
at the present.%
\footnote{We thank C.G.~Tamm (IMPA) for enlightening discussions on
this combinatorial exercise.}
The treatment of the Neumann boundary conditions and its probabilistic
interpretation goes beyond the scope of \cite{Gr92} albeit it shows the
power of ideas presented.

%-----------------------------------------------------------------------------
\section{Numerical approach} \label{sec:numeric}

In this section we consider a numerical approach based on level set methods
for the inverse doping profile problem in the stationary linearized unipolar
case close to equilibrium (see Identification Problem~\ref{probl:unipol}).
We compare our results with the ones obtained in \cite{BELM04}, where a
Landweber-Kaczmarz iterative method was used to reconstruct the doping
profile function.

\subsection*{Framework}

As already mentioned in Section~\ref{sec:model}, the main task in this
inverse problem consists in the identification of the coefficient $\gamma$
in the elliptic PDE
\begin{equation} \label{eq:model-DtN}
{\rm div} (\gamma \nabla u) \, = \, 0 \ \ {\rm in }\ \Omega\, ,\ \ \ \ \ 
   u \, = \, U \ \ {\rm on }\ \partial\Omega_D\, ,\ \ \ \ \ 
   u_\nu \, = \, 0 \ \ {\rm on }\ \partial\Omega_N \, .
\end{equation}
Thus, we can reduce the inverse doping profile problem to the problem of
identifying a piecewise constant function $\gamma(x)$ in (\ref{eq:model-DtN})
from measurements of the DtN map
\begin{equation} \label{eq:DtN-map}
\Lambda_\gamma : \begin{array}[t]{rcl}
  H^{3/2}(\partial\Omega_D) & \to & H^{1/2}(\Gamma_1) \\
  U & \mapsto & \gamma\, u_\nu |_{\Gamma_1}
\end{array}
\end{equation}
(for simplicity we assume $\gamma(x) \in \{ 1, 2 \}$ a.e. in $\Omega$).

Notice that, due to the nature of the boundary conditions related to
practical experiments, we have to restrict the domain of definition
of the DtN operator to the linear subspace $ D(\Lambda_\gamma) :=
\{ U \in H^{3/2}(\partial\Omega_D); \, U|_{\Gamma_1} = 0\}$. Furthermore,
the measurements (Neumann data) are only available at $\Gamma_1$.
This is the essential difference between the parameter identification
problem in (\ref{eq:model-DtN}) and the classical inverse problem in {\em
electrical impedance tomography}, namely the fact that both Dirichlet and
Neumann data are known only at specific parts of the boundary.

For this particular type of DtN operators there are so far no analytical
results concerning identifiability and, to our knowledge, the few numerical
results in the literature are those discussed in \cite{BELM04, BEMP01,
BEM02, FIR02}. 

In this section we shall work within the following framework:
\begin{itemize}
\item Parameter space: \ $\mathcal X := L^2(\Omega)$;
\item Input (fixed): \ $U_j \in H^{3/2}(\partial\Omega_D)$, \, with \,
$U_j |_{\Gamma_1} = 0$, \ $1 \le j \le N$;
\item Output (data): \ $Y = \big\{ \Lambda_\gamma(U_j) \big\}_{j=1}^N
\in \big[L^2(\Gamma_1)\big]^N =: \mathcal Y$;
\item Parameter to output map: \ 
$F: \begin{array}[t]{rcl}
      D(F) \subset \mathcal X & \to & \mathcal Y \\
      \gamma(x) & \mapsto & \big\{ \Lambda_\gamma(U_j) \big\}_{j=1}^N
    \end{array}$
\end{itemize}
where the domain of definition of the operator $F$ is
$$ D(F) := \{ \gamma \in L^2(\Omega) ; \, \gamma_+ \ge \gamma(x)
           \ge \gamma_- > 0, \mbox{ a.e. in } \Omega \} $$
(here $\gamma_-$ and $\gamma_+$ are appropriate constants).
We shall denote noisy data by $Y^\delta$ and assume that the data error
is bounded by $\| Y - Y^\delta \| \le \delta$.
Thus, we are able to represent the inverse doping problem in the general form
\begin{equation} \label{eq:ip-dd}
F(\gamma) \, = \, Y^\delta .
\end{equation}

For the concrete numerical test performed in this section as well as in
\cite{BELM04}, $\Omega \subset \mathbb R^2$ is the unit square and the
boundary parts are
\begin{equation*}
\Gamma_1 \, := \,  \{ (x,1) \, ;\ x \in (0,1) \} \, ,\ \ \
\Gamma_0 \, := \, \Gamma_1 \cup \{ (x,0) \, ;\ x \in (0,1) \} \, ,\ \ \
\partial\Omega_D \, := \, \Gamma_0 \cup \Gamma_1\, ,
\end{equation*}
\begin{equation*}
\partial\Omega_N \ := \ \{ (0,y) \, ;\ y \in (0,1) \} \cup
   \{ (1,y) \, ;\ y \in (0,1) \} \, .
\end{equation*}
The fixed inputs $U_j$ vanish at $\Gamma_1$ and are chosen to be piecewise
constant functions on $\Gamma_0 = \{ (x,0) \, ;\ x \in (0,1) \}$.
\begin{equation*}
U_j(x) \ := \ \left\{ \begin{array}{rl}
1, & |x - x_j| \le \delta x \\
0, & {\rm else} \end{array} \right.
\end{equation*}
where the points $(x_j,1)$, $j=1,\dots,N$, are uniformly distributed in
the segment $\Gamma_0$.

The next lemma describes some crucial properties of the operator $F$, that
will be necessary for the analysis of the iterative methods discussed in
this paper. Here, only a sketch of the proof of Lemma~\ref{lem:pto-wd} is
given, for details see \cite{BEMP01}.

\begin{lemma} \label{lem:pto-wd}
Let the voltages $\{ U_j \}_{j=1}^N$ be chosen in the neighborhood of
$U = 0$. The parameter-to-output map $F$ defined above is well-defined
and Fr\'echet differentiable on $D(F)$.

{\it Proof:} \rm
The first statement follows from the well-definedness of the V-C map, cf.
Lemma~\ref{lem:vc-wd}. The Fr\'echet differentiability of $F$ follows from
the differentiability of the V-C map (see Lemma~\ref{lem:vc-wd}) together
with the differentiability of the map that takes the doping profile $C$
onto the solution $(V,u,v)$ of (\ref{eq:dd-sys}). \hfill $\square$
\end{lemma}

\subsection*{A competing approach: Landweber-Kaczmarz method}

In \cite{BELM04} a Landweber-Kaczmarz method was used to reconstruct the
doping profile function. This corresponds to an iterative method of
steepest descent type for solving the least square formulation of the
inverse problem.

A simple and robust iterative method to solve the inverse problem
(\ref{eq:ip-dd}) is the {\em Landweber iteration} \cite{DES98, EHN96, ES00,
HNS95}. This iteration is known to generate a regularization method for the
inverse problem, the stopping index playing the role of the regularization
parameter (for regularization theory see, e.g., \cite{EHN96, EKN89, ES00,
TA77}). 

The {\em Landweber-Kaczmarz method} \cite{KS02} results from the
coupling of the strategies of both the Landweber and the Kaczmarz
iterative methods.
The  motivation for this choice of strategy lays in the fact that the
data in (\ref{eq:ip-dd}) consists of a vector of measurements
$\{ \Lambda_\gamma(U_j) \}_{j=1}^N$ and the principal characteristic
of the Kaczmarz method is the minimization, at each iteration step, of
a least square functional that takes into account only one component
of this measurement vector.
It is worth mentioning that this method has already been successfully
applied to {\em electrical impedance tomography} by Nachman in \cite{Nt86}.

To formulate the method, we first need to define the components of the
parameter to output map: $F = \{ \mathcal F_j\}_{j=1}^N$, where
$\mathcal F_j: L^2(\Omega) \supset D(F) \ni \gamma \mapsto
\Lambda_\gamma(U_j) \in L^2(\Gamma_1)$.
Now, setting $Y_j^\delta := \mathcal F_j(\gamma^\delta)$,
$1 \le j \le N$, the Landweber-Kaczmarz iteration can be written in the form
\begin{equation} \label{eq:land-kacz}
  \gamma^\delta_{k+1} \, = \, \gamma_k^\delta
                          - \mathcal F_k'(\gamma_k^\delta)^*
  \big( \mathcal F_k (\gamma_k^\delta) - Y^\delta_k \big)\, ,
\end{equation}
for $k = 1, 2, \dots$, where we adopted the notation
$\mathcal F_k := \mathcal F_j, \ \ Y^\delta_k := Y^\delta_j$,
whenever $k = i \, N + j$, and $i = 0, 1, \dots$, and $j = 1, \dots, N$.

Notice that each step of the Landweber-Kaczmarz method consists of one
Land\-weber iterative step with respect to the $j$-th component of the
residual in (\ref{eq:ip-dd}). These Landweber steps are performed in a
cyclic way, using the components of the residual $\mathcal F_j(\gamma)
- Y^\delta_j$, $1 \le j \le N$, one at a time.

\subsection*{A level set approach}

In this paper we propose a level set type method to approximate the
solution of (\ref{eq:model-DtN}).
In the sequel, the function spaces $\mathcal X$, $\mathcal Y$ as well as
the operators $F$, $\Lambda_\gamma$ and also the sets $\Omega$, $\partial
\Omega_D$, $\partial \Omega_N$, $\Gamma_0$, $\Gamma_1$ are the same as
above. We assume, however, that only one measurement is given, i.e., only
one pair of voltage-current data is available for the reconstruction.
This assumption corresponds to the choice $N=1$ in the definition of the
space $\mathcal Y$.

Our numerical approach is based on the level set method introduced in
\cite{LS03,FSL05}. According to this strategy, one represents the unknown
P-N junction by the zero level set of an $H^1$-function $\phi: \Omega \to
\mathbb R$, in such a way that $\phi(x) > 0$ if $\gamma(x)=2$ and
$\phi(x) < 0$ if $\gamma(x)=1$.
Starting from some initial guess $\phi_0 \in H^1(\Omega)$, one
solves the Hamilton-Jacobi equation
\begin{equation} \label{eq:ls-hj}
\frac{\partial \phi}{\partial t} + V \nabla \phi = 0
\end{equation}
where $V = - v \, \nabla\phi / |\nabla\phi|^2$ and the {\em velocity}
$v$ solves
\begin{equation} \label{eq:ls-vel}
\left\{ \!\! \begin{array}{l}
  (\Delta - I) v = \frac{\delta(\phi(t))}{|\nabla \phi(t)|}
  \left[F'(\chi(t))^*(F(\chi(t)) - Y^\delta) - \beta \nabla \!\!\cdot\!\!
  \left( \frac{\nabla P(\phi)}{|\nabla P (\phi)|} \right) \right]
  ,\, {\rm in} \ \Omega \\
  \D \frac{\partial v}{\partial\nu} = 0 \ ,\, {\rm on}\ \partial\Omega \, .
\end{array} \right.
\end{equation}
Here, $\alpha>0$ is a regularization parameter and $\chi = \chi(x,t)$ is
the projection of the level set function $\phi(x,t)$ defined by:
\begin{equation*}
\chi(x,t) = P(\phi(x,t)) :=
            \left\{ \begin{array}{ll}
              2, & {\rm if}\ \phi(x,t) > 0 \\
              1, & {\rm if}\ \phi(x,t) < 0
            \end{array} \right. .
\end{equation*}

In \cite{FSL05,LS03} this level set method was used to reconstruct inclusions
$D \subset\subset \Omega$. Notice that, in our case, the set $D$ corresponds
to the P-region (see Figure~\ref{fig:diode}) and the condition $\overline{D}
\subset \Omega$ is not satisfied. This fact, however, does not affect the
derivation of the Hamilton-Jacobi equation (\ref{eq:ls-hj}). Moreover, it
does not affect the derivation of the boundary conditions for the elliptic
problem (\ref{eq:ls-vel}) either.

The family $\chi(\cdot,t)$ approximate the doping profile $\gamma(\cdot)$ as
$t \to \infty$. This follows from the fact that the solution $\phi(\cdot,t)$
of (\ref{eq:ls-hj}) converges to the minimum of the Tikhonov functional
\begin{equation} \label{eq:Ga}
\Ga(\phi) := \|F(P(\phi)) - Y^\delta\|^2_{\cal Y} +
   \alpha \big( 2\beta |P(\phi)|_\bv + \|\phi - \phi_0\|^2_{H^1(\Omega)} \big)
\end{equation}
as $t \to \infty$, for each regularization parameter $\alpha>0$ ($\beta>0$ is
fixed). See \cite[Definition~2.2]{FSL05} for the precise definition of a
minimizer of $\Ga(\phi)$.

The next lemma corresponds to specific results selected from \cite{FSL05}.
It allows a better understanding of the least-square problem behind the
level set formulation and also analytically substantiates the numerical
results presented in the sequel.

\begin{lemma} \label{lem:ls-wp}
Stability, Convergence and Well-Posedness:
\begin{itemize}
\item[{\bf (a)}] Let $Y^\delta = Y$ (noiseless case), and let $\phi_\alpha$
be a minimizer of $\Ga$.
Then, for every sequence $\{\alpha_k\}_{k \in \mathbb N}$ converging to zero,
there exists a subsequence $\{ \alpha_{k(l)} \}_{l \in \mathbb N}$, such that
$\{ \phi_{ \alpha_{k(l)} } \}_{l \in \mathbb N}$ is strongly convergent.
Moreover, the limit is a minimal norm solution of (\ref{eq:ip-dd}).
\item[{\bf (b)}] Let $\|Y^\delta - Y\|_{\mathcal Y} \leq \delta$. If
$\alpha = \alpha(\delta)$ satisfies $\D\lim_{\delta \to 0} \alpha(\delta)
= 0$ and $\D\lim_{\delta \to 0} \T\frac{\delta^2}{\alpha (\delta)} = 0$,
then, for a sequence $\{\delta_k\}_{k\in\mathbb N}$ converging to 0,
the sequence $\phi_{\alpha(\delta_k)}$ converges to a minimal norm
solution of (\ref{eq:ip-dd}).
\item[{\bf (c)}] For any given $\phi_0 \in H^1(\Omega)$ the functional
$\Ga$ attains a minimizer.
\end{itemize}
\end{lemma}

\begin{remark}[Level set algorithm]
For the reader's convenience, we briefly describe the level set algorithm
related to (\ref{eq:ls-hj}), (\ref{eq:ls-vel}). Here, $P_\ve$ is the
approximation defined in \cite[Section~2]{FSL05} for the operator $P$.
The adjoint operator $(F')^*$ as well as its evaluation on a given
vector is derived in \cite[Section~4]{BELM04}.

\begin{tt} \begin{enumerate}
\item[{\bf 1.}] \ Evaluate the residual \
$r_k := F( P_\ve(\phi_k) ) - Y^\delta$;
\item[{\bf 2.}] \ Evaluate\, $w_k := F'( P_\ve(\phi_k) )^* (r_k)$;
\item[{\bf 3.}] \ Evaluate $v_k \in H^1(\Omega)$, satisfying
$$ \begin{aligned}
(\Delta - I) v_k &= P_\ve'(\phi_k)
   \left( w_k - \beta P_\ve'(\phi_k) \nabla \!\!\cdot\!\!
          \left( \T\frac{\nabla P_\ve(\phi_k)}{|\nabla P_\ve(\phi_k)|} \right)
   \right) \, , {\tt in}\ \Omega \\
\partial v_k / \partial \nu & = 0 \, , {\tt on}\ \partial\Omega \; ;
\end{aligned} $$
\item[{\bf 4.}] \ Update the level set function \ $\phi_{k+1} =
\phi_k + v_k$.
\end{enumerate} \end{tt}
\end{remark}

\begin{figure}[b]
\centerline{ \epsfysize3.4cm \epsfbox{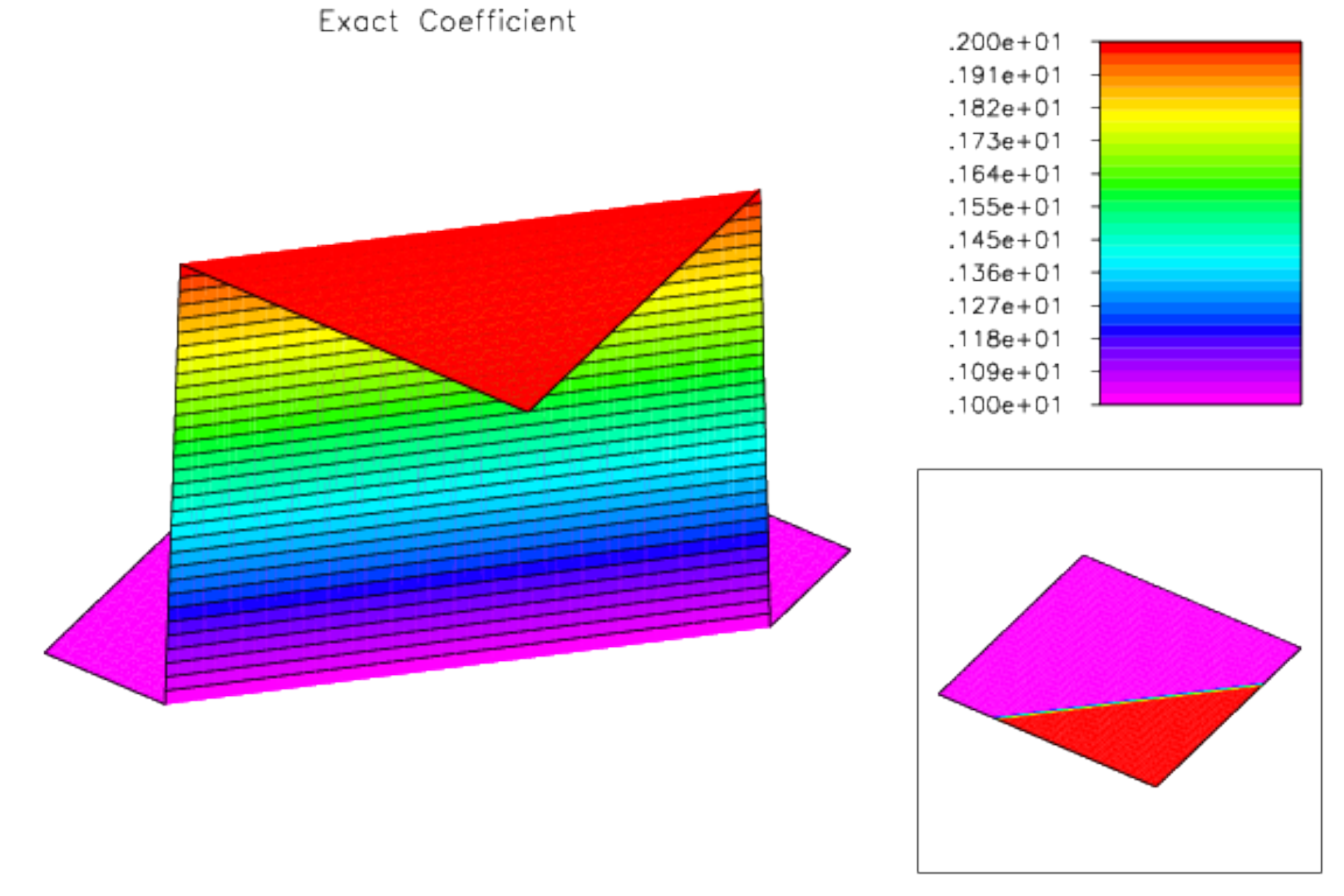} \hfil
             \epsfysize3.4cm \epsfbox{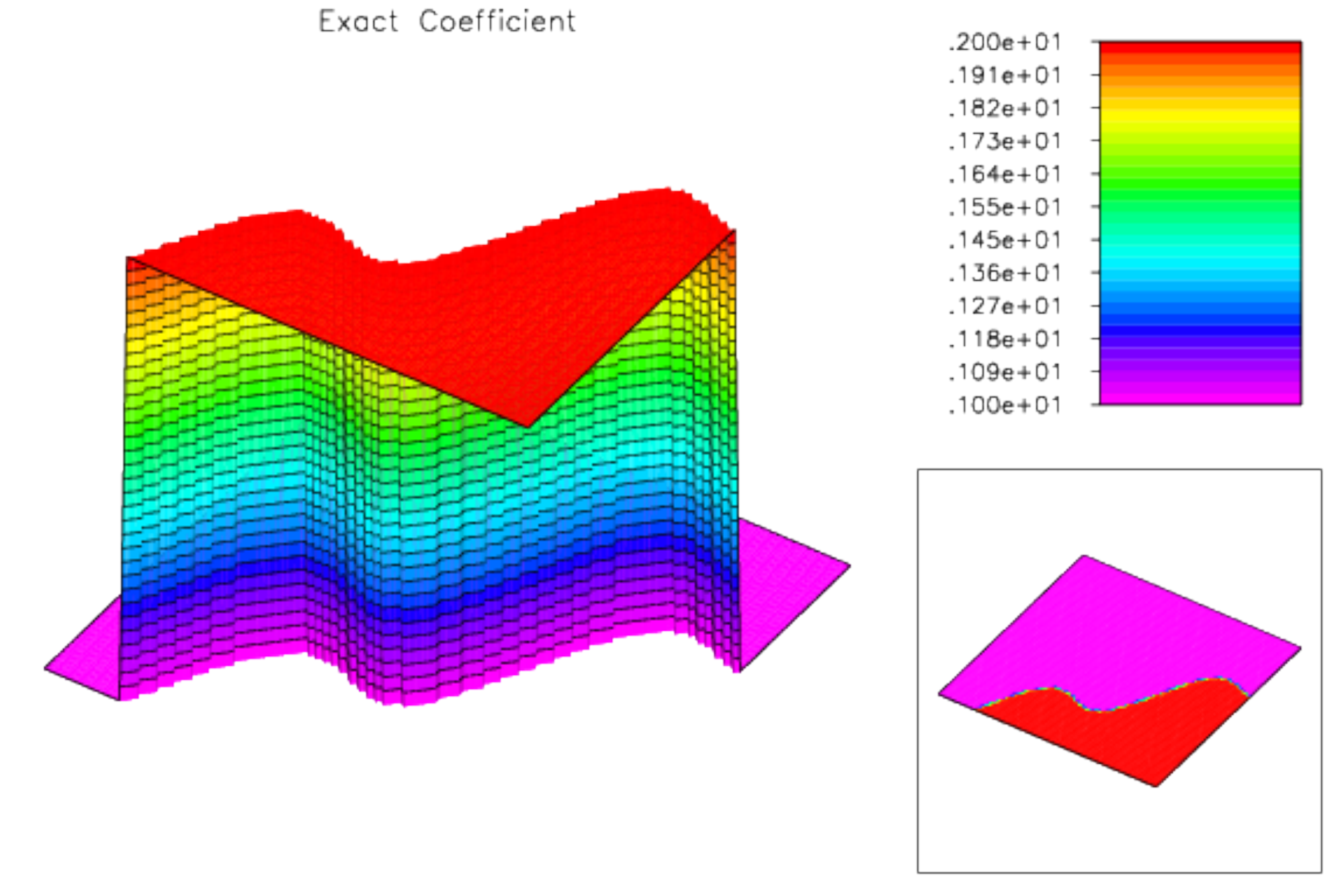} \hfil
             \epsfysize3.4cm \epsfbox{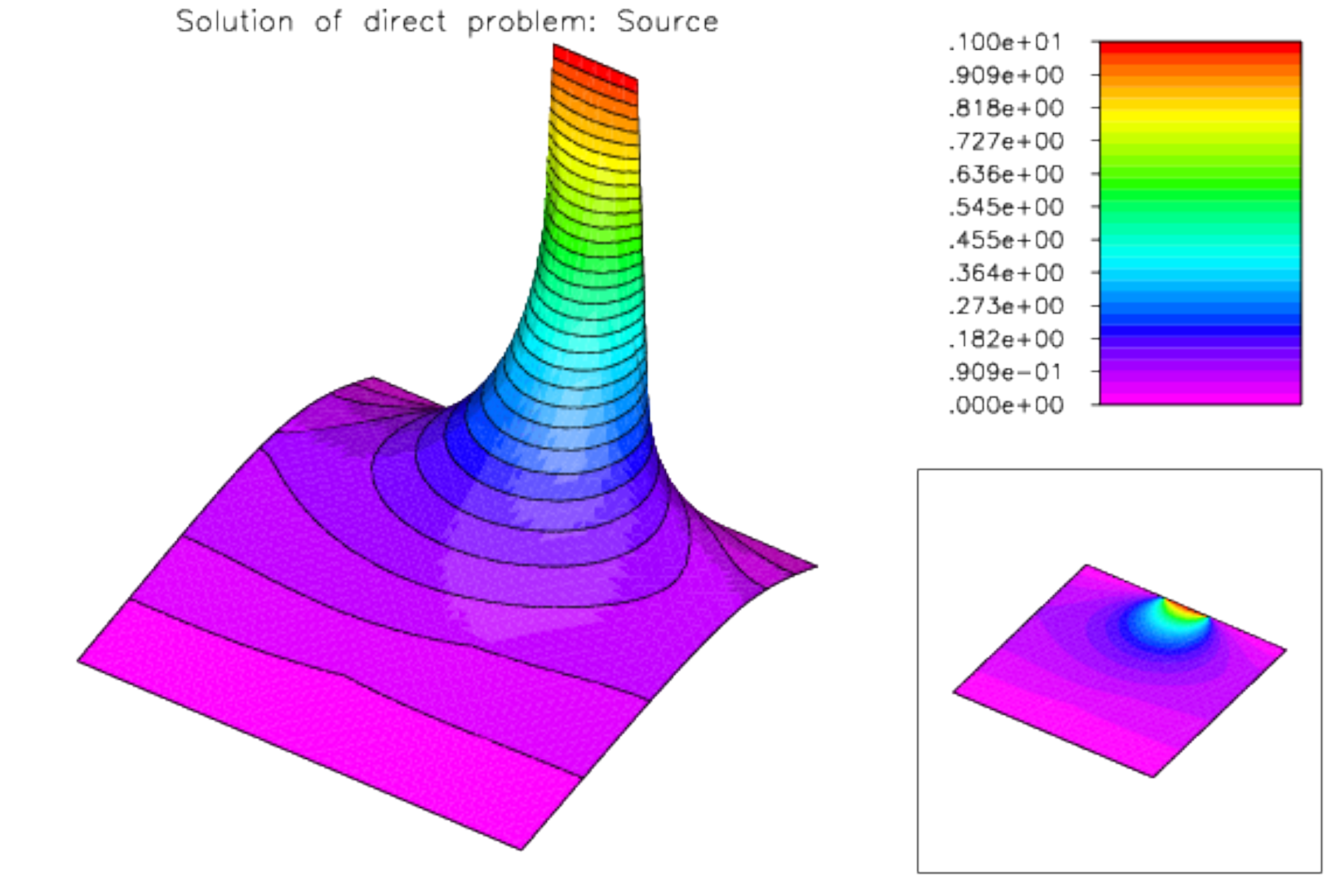}  }
\centerline{\hfil (a) \hskip3.8cm  (b) \hskip3.8cm  (c) \hfil}
\caption{\small Pictures (a), (b) show the doping profiles to be
reconstructed in the two different experiments for the level set method.
On picture (c) the problem data is shown: The source $U(x)$ appears as
the Dirichlet boundary condition at $y=1$ ($\Gamma_0$ is the upper right
edge). The corresponding current is measured at the contact $\Gamma_1$
(lower left edge), where $U(x)$ is assumed to vanish.}
\label{fig:ls-setup}
\end{figure}

\begin{figure}[t]
\centerline{ \epsfysize3.4cm \epsfbox{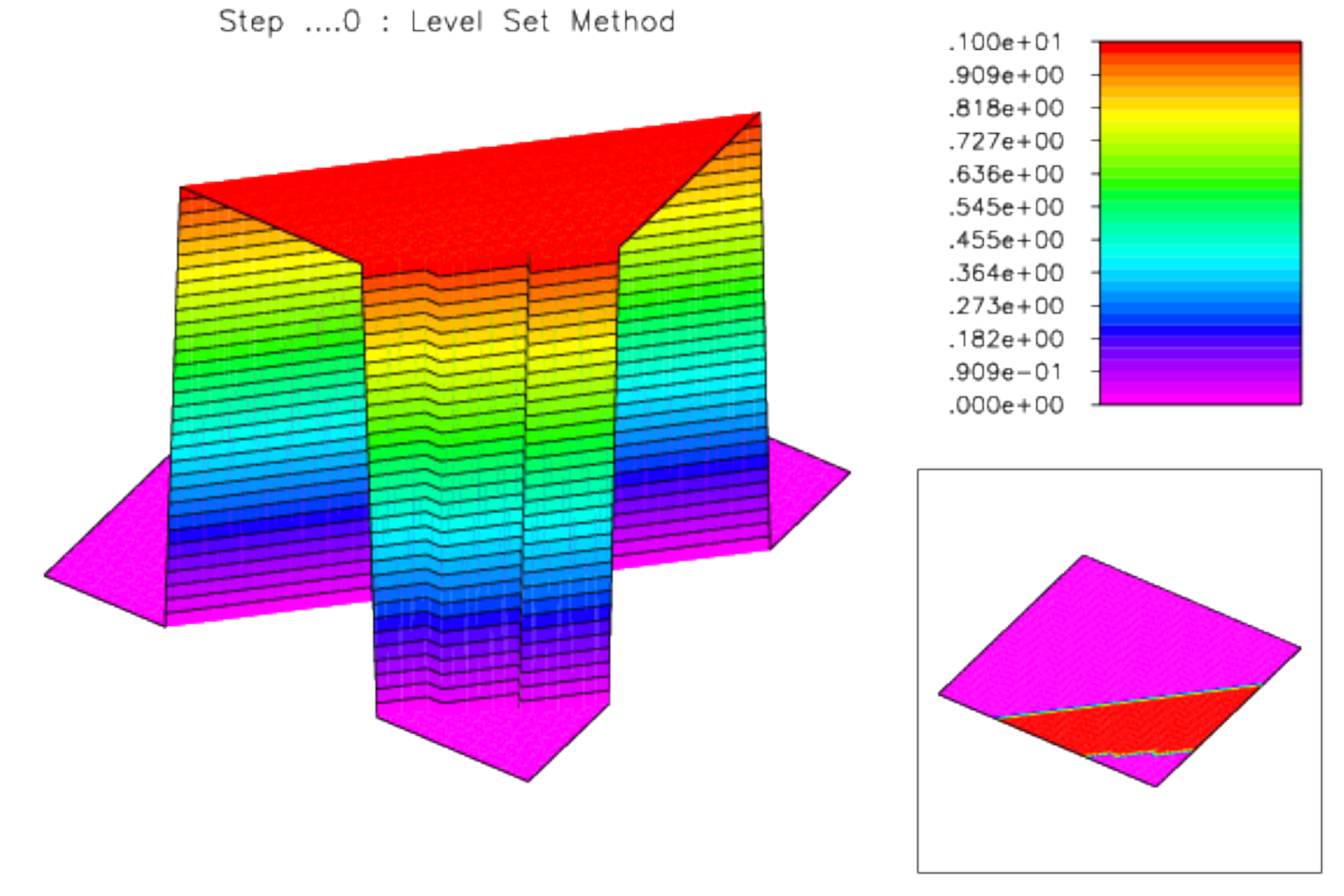} \hfil
             \epsfysize3.4cm \epsfbox{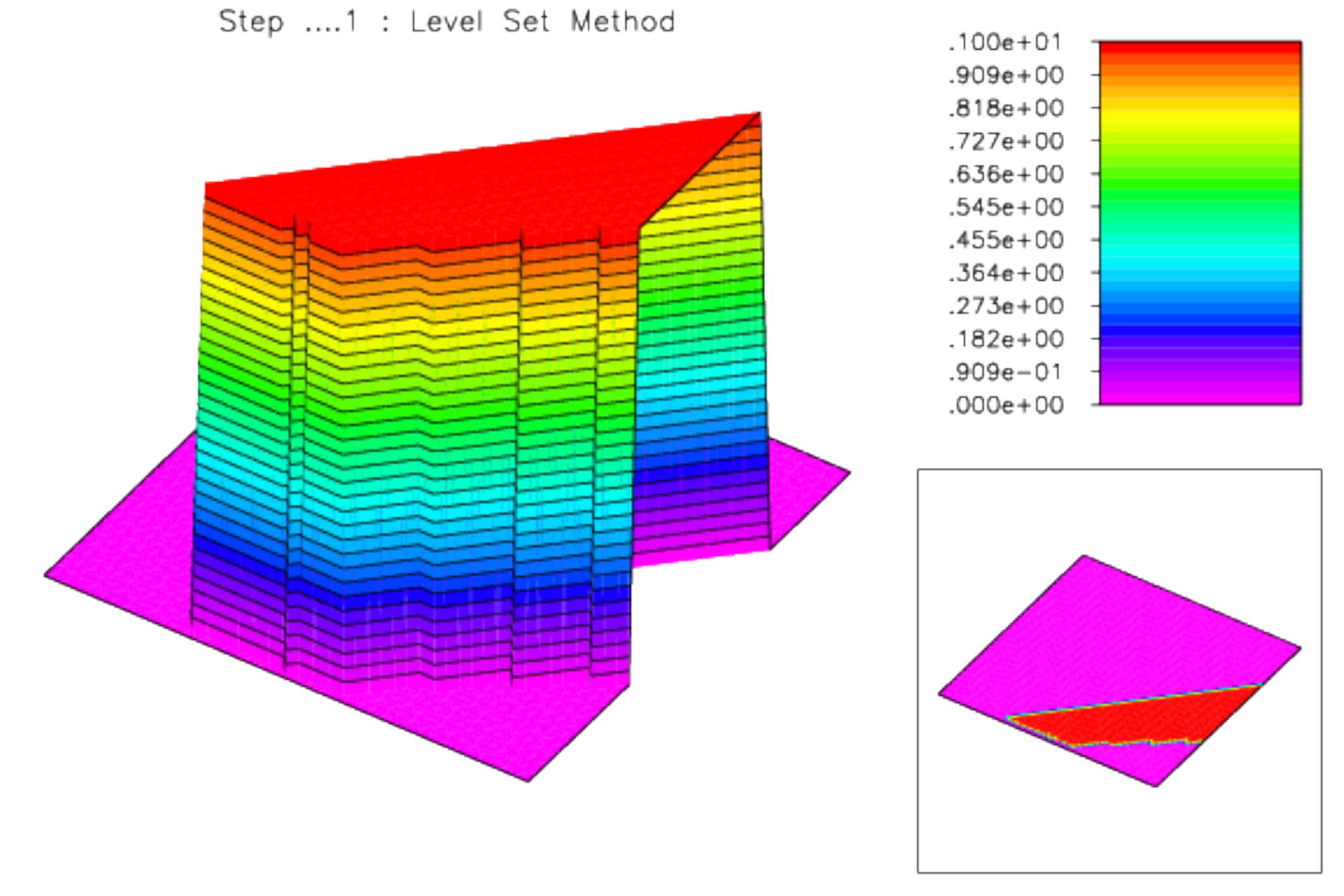} \hfil
             \epsfysize3.4cm \epsfbox{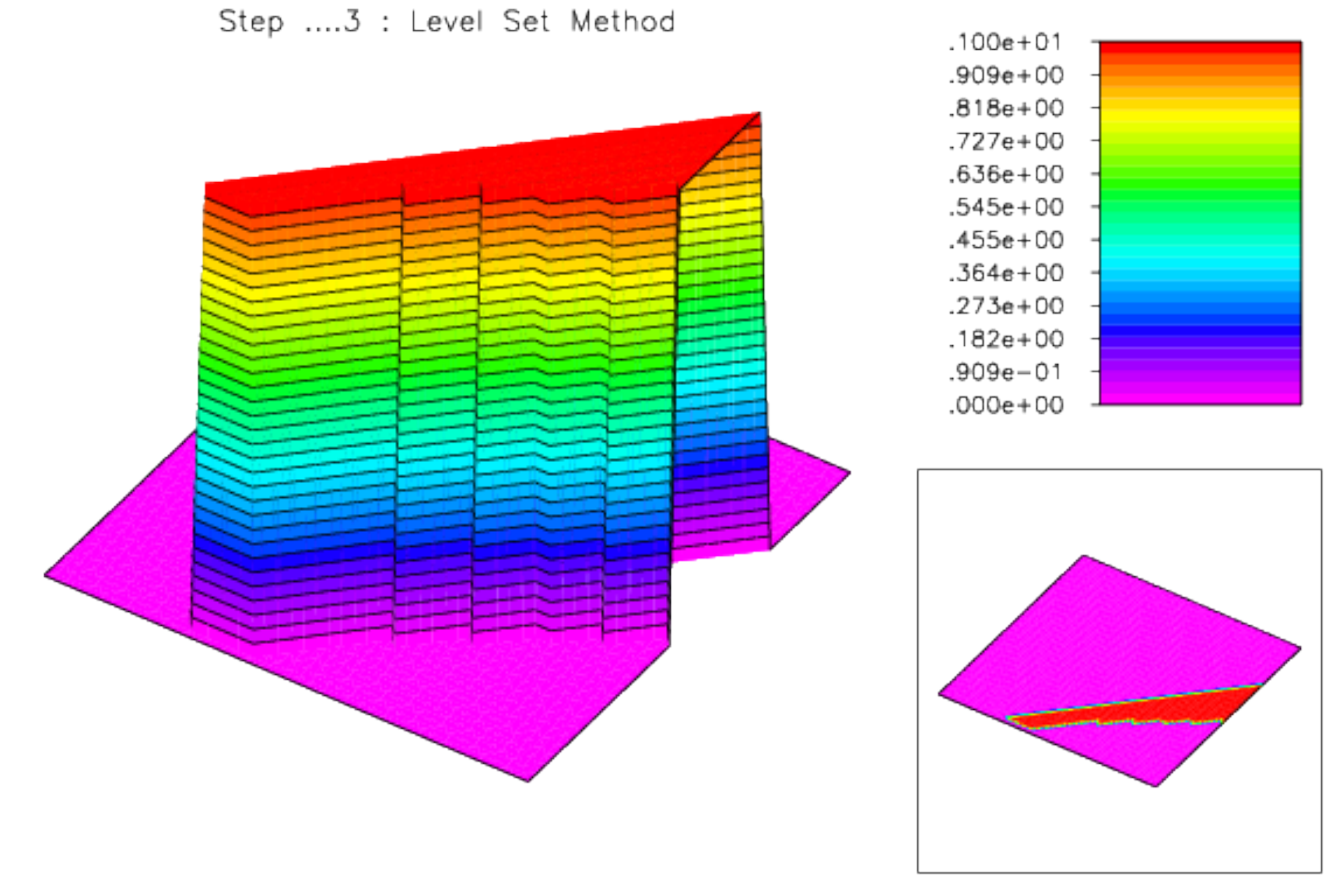} }
\bigskip
\centerline{ \epsfysize3.4cm \epsfbox{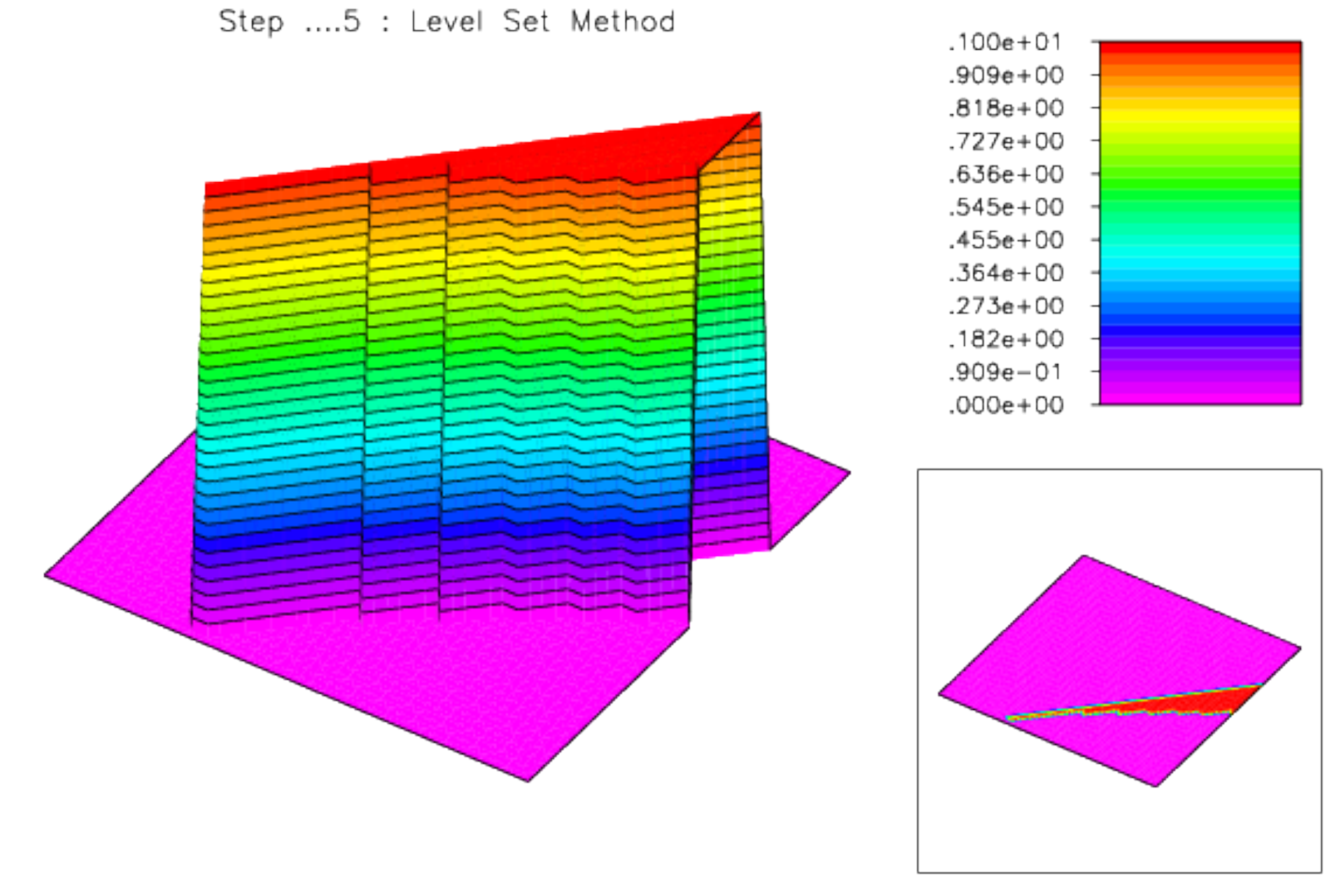} \hfil
             \epsfysize3.4cm \epsfbox{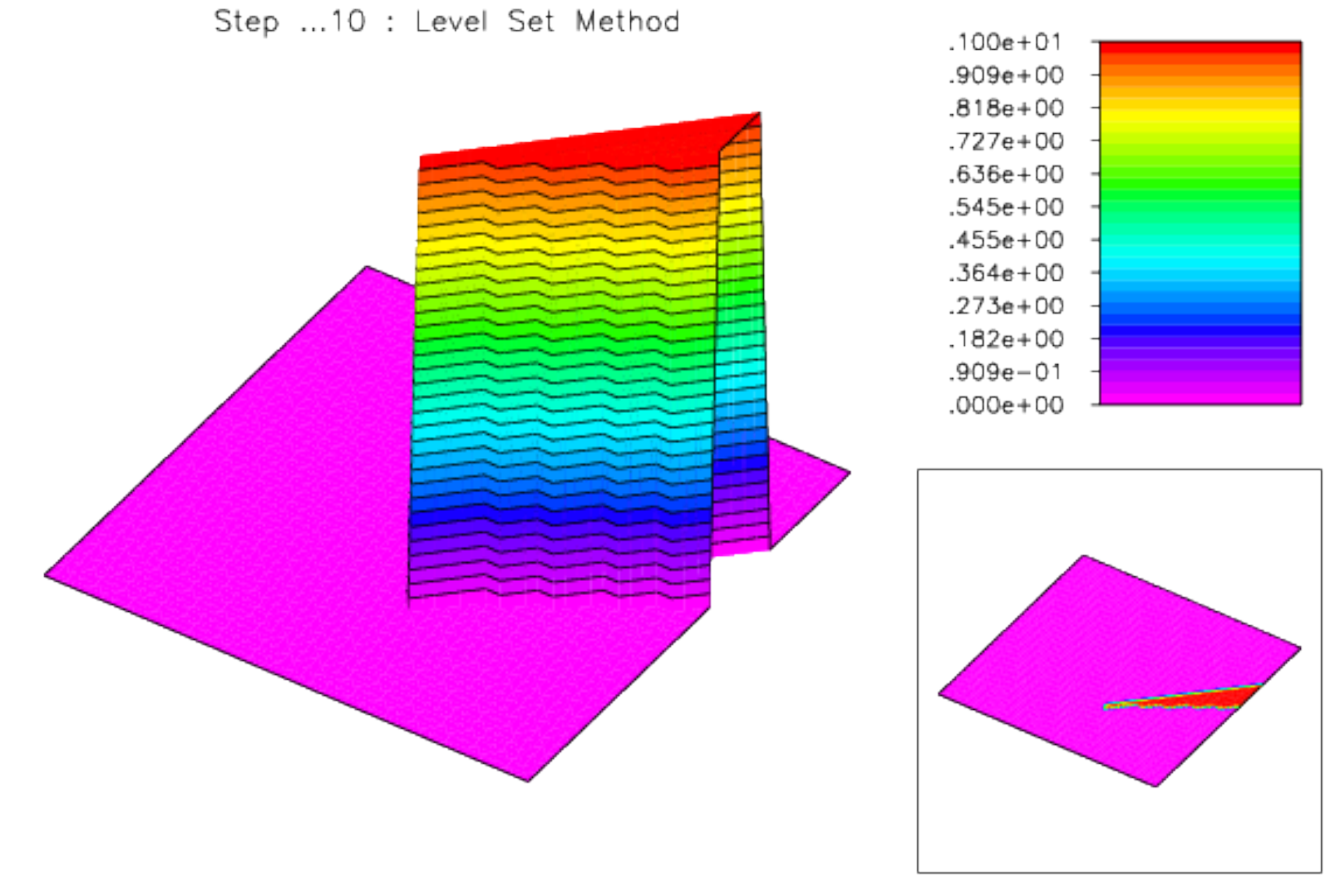} \hfil
             \epsfysize3.4cm \epsfbox{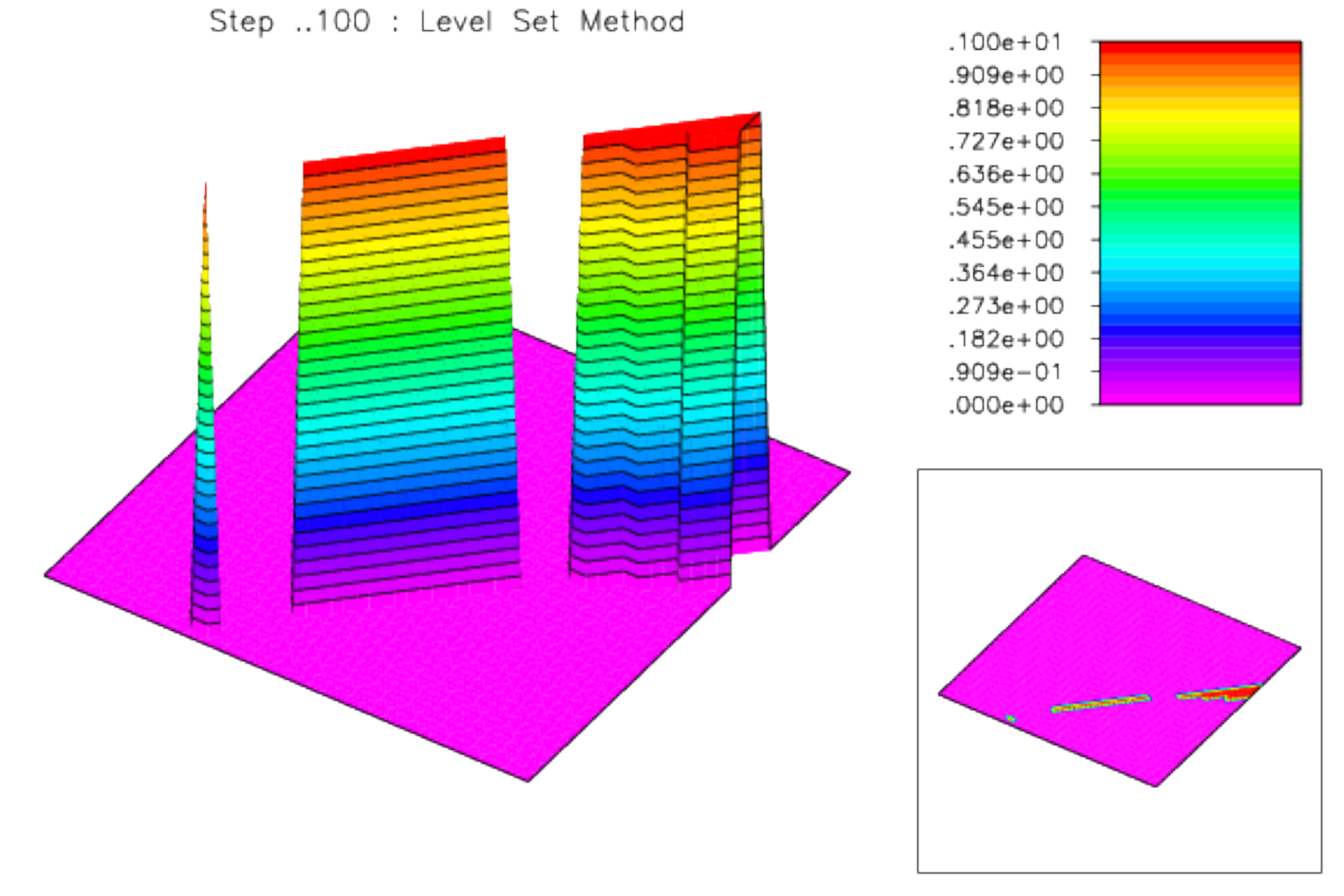} }
\caption{\small First numerical experiment (linear P-N junction): Evolution
of the iteration error for the level set method and exact data.}
\label{fig:evol-ls1}
\end{figure}

\begin{figure}[ht]
\centerline{ \epsfysize3.4cm \epsfbox{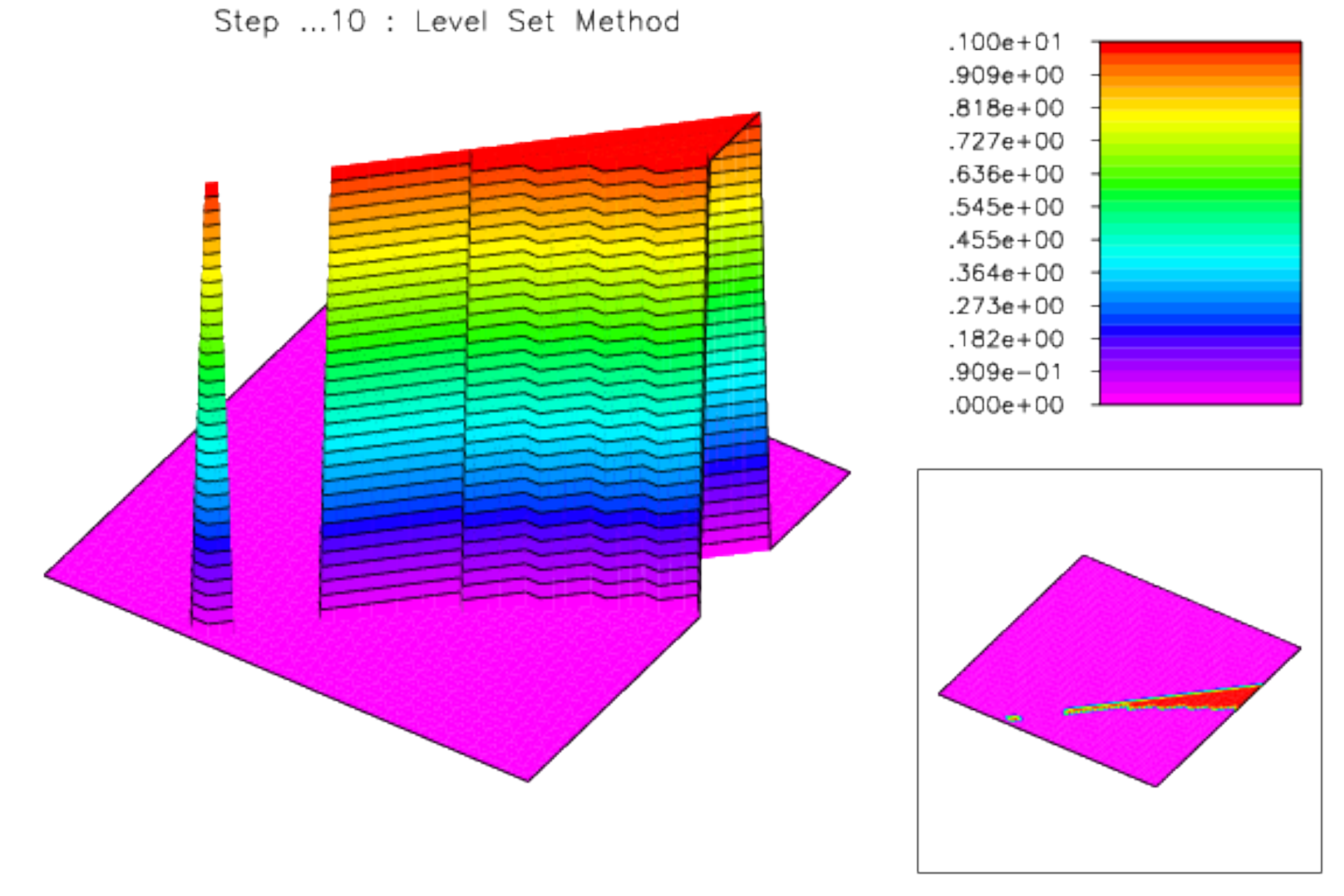} \hfil
             \epsfysize3.4cm \epsfbox{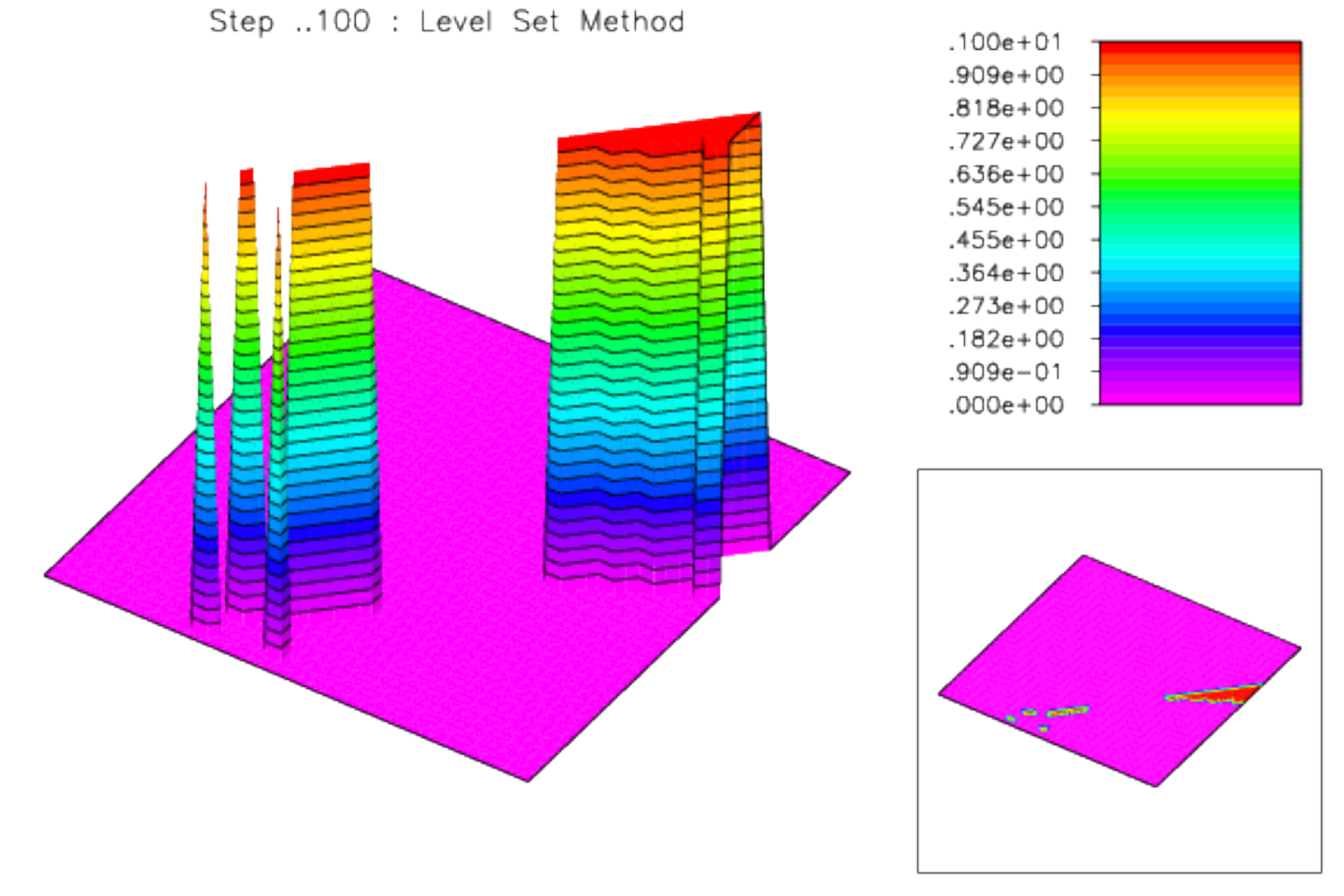} \hfil
             \epsfysize3.4cm \epsfbox{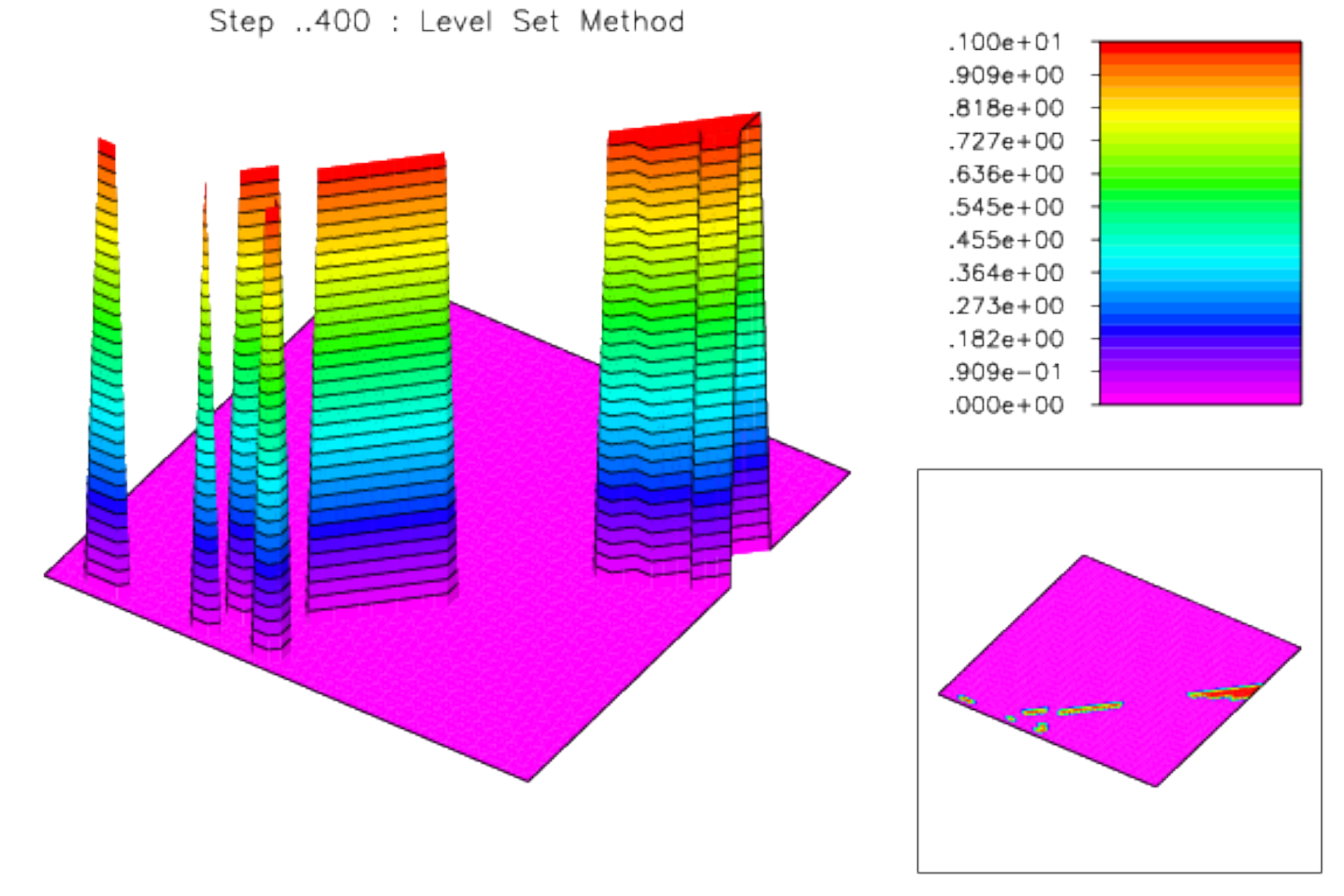} }
\caption{\small First numerical experiment (linear P-N junction): Evolution of
the iteration error for the level set method and data with 10\% random noise.}
\label{fig:evol-ls2}
\end{figure}

\begin{figure}[ht]
\centerline{ \epsfysize3.4cm \epsfbox{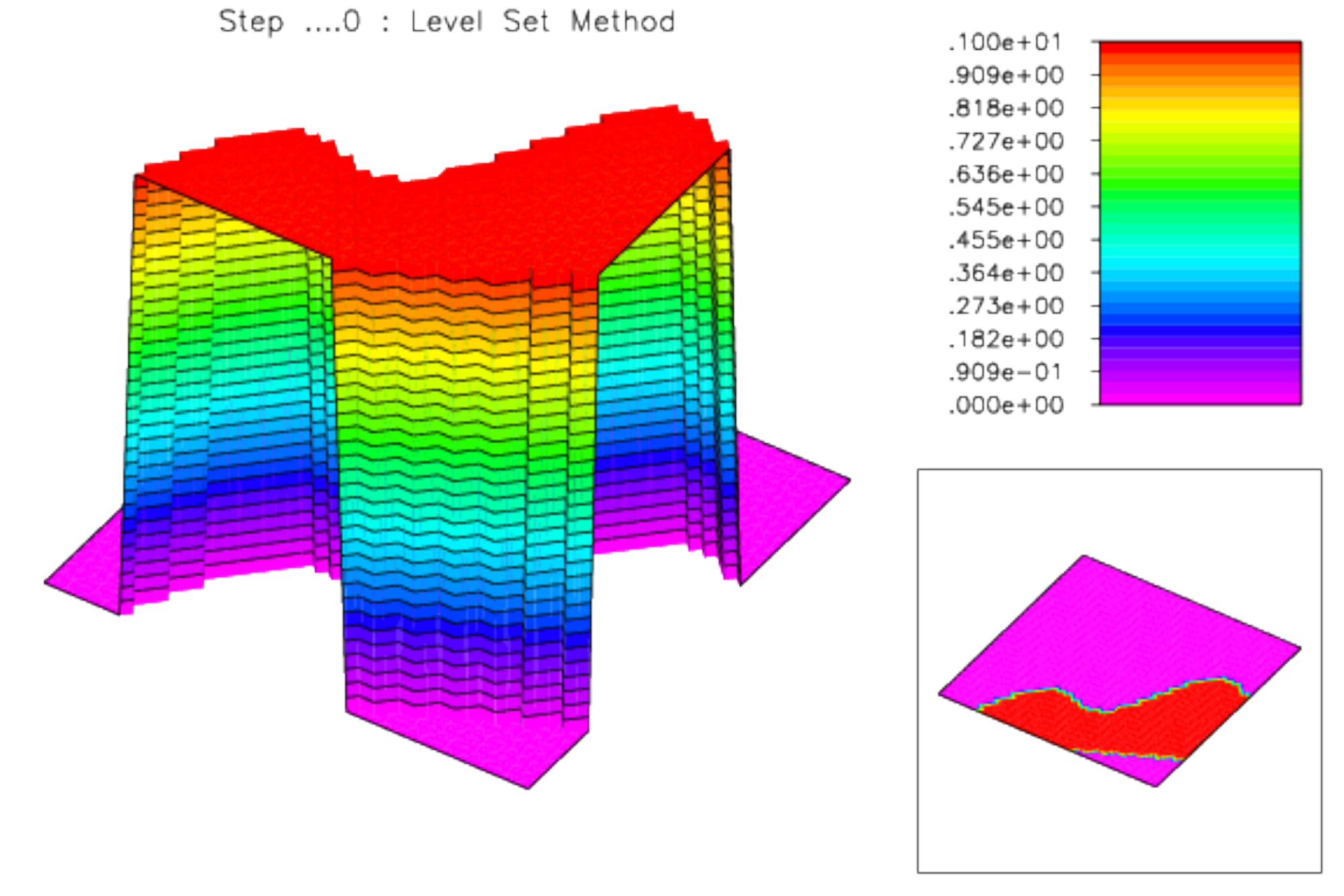} \hfil
             \epsfysize3.4cm \epsfbox{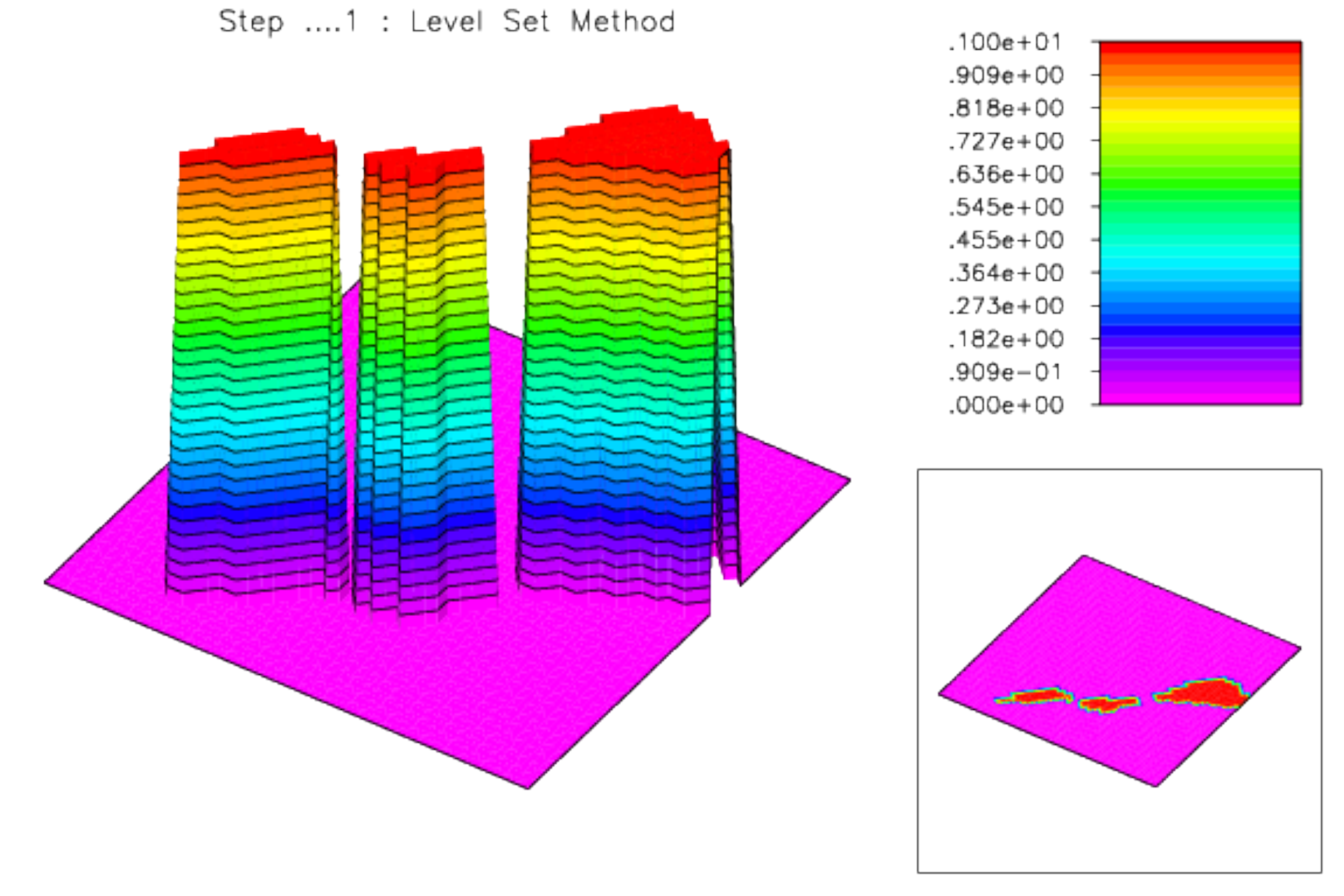} \hfil
             \epsfysize3.4cm \epsfbox{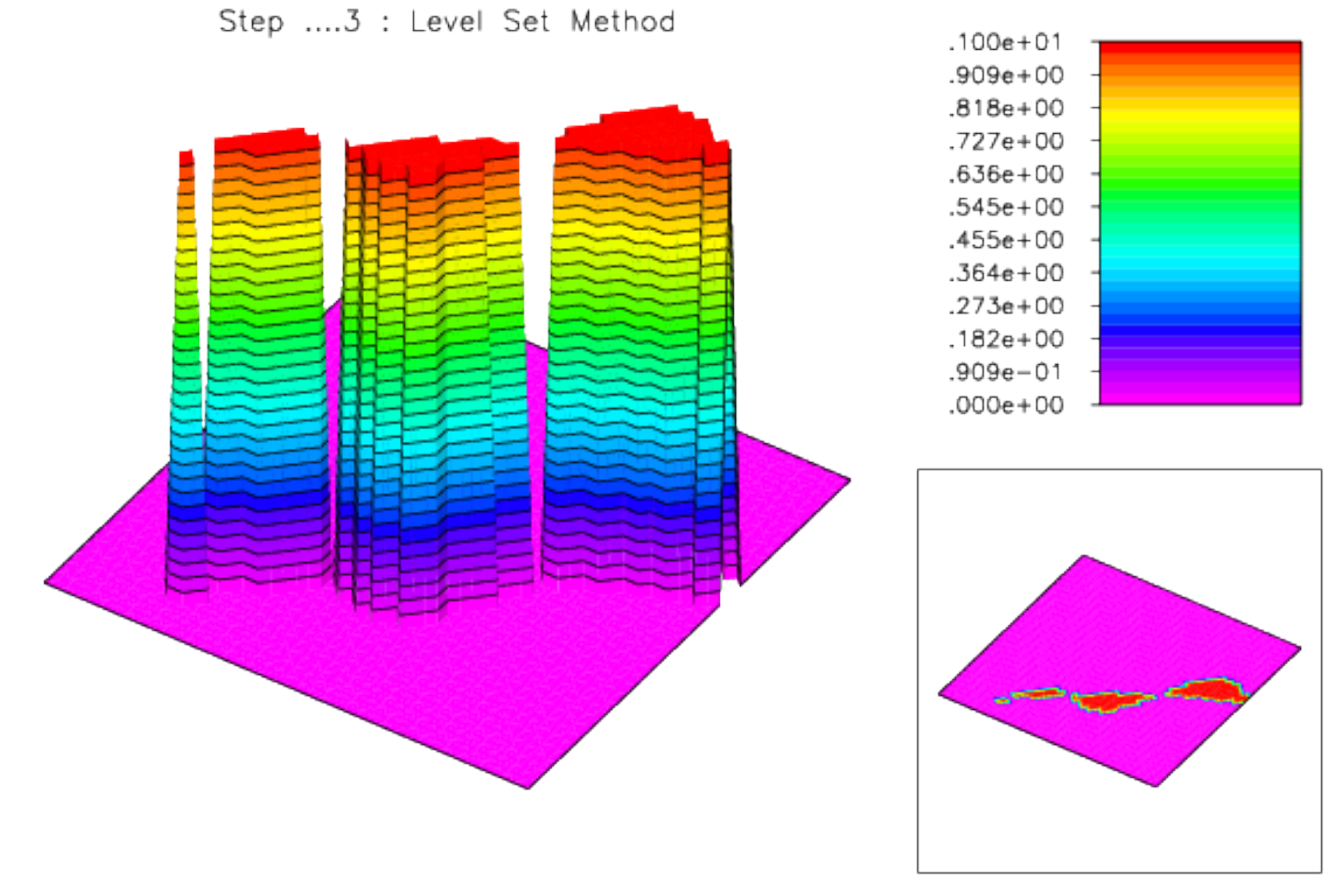} }
\bigskip
\centerline{ \epsfysize3.4cm \epsfbox{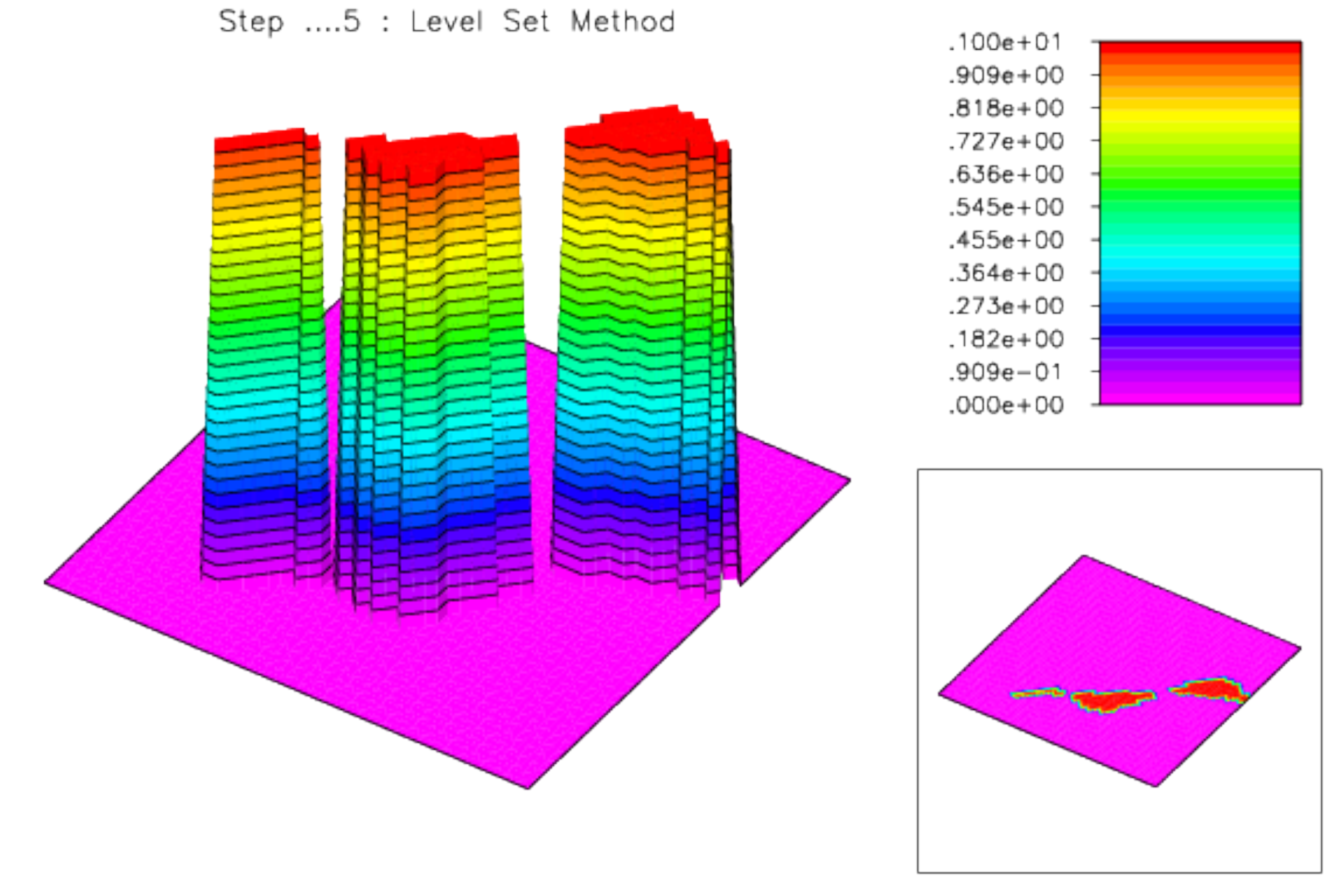} \hfil
             \epsfysize3.4cm \epsfbox{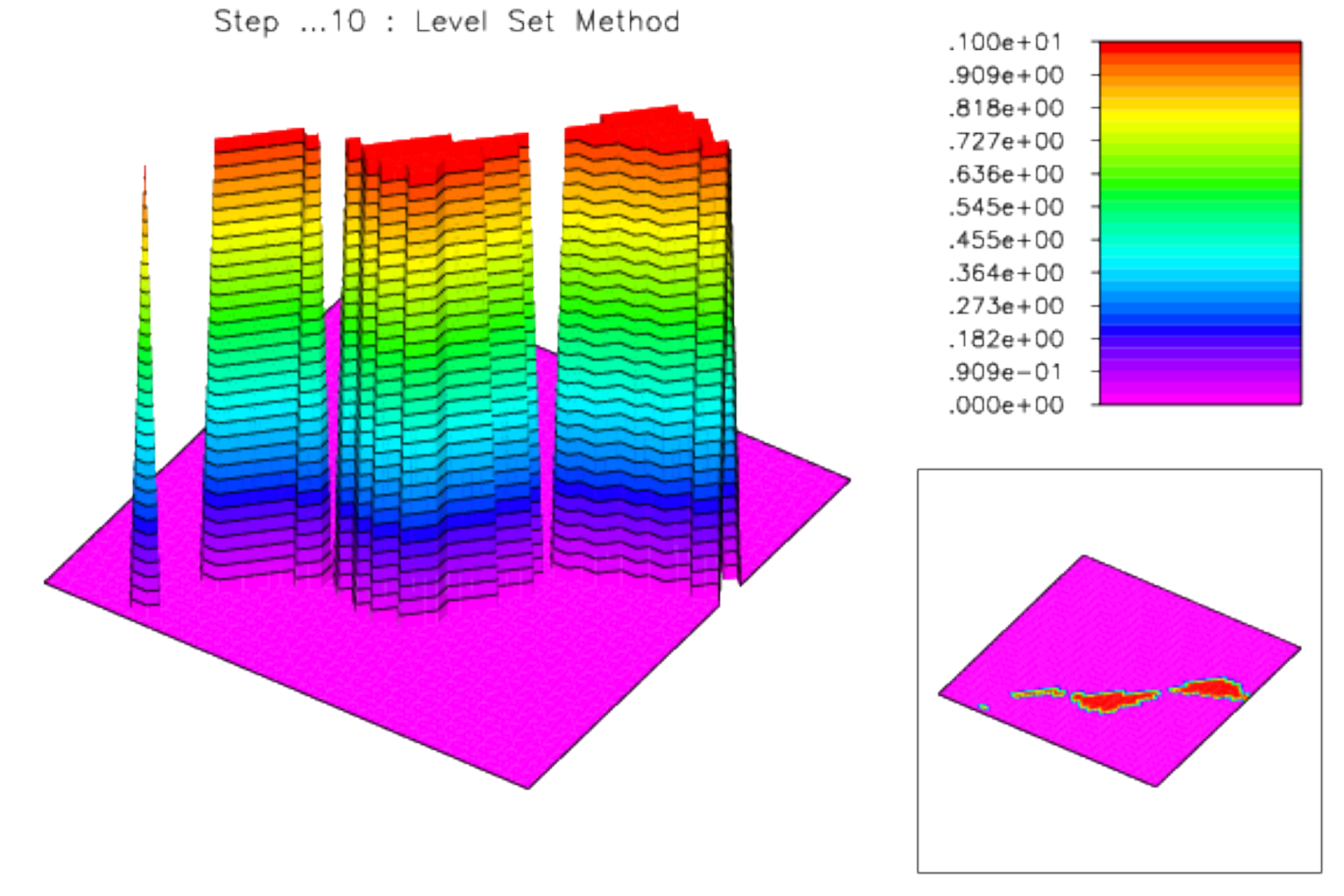} \hfil
             \epsfysize3.4cm \epsfbox{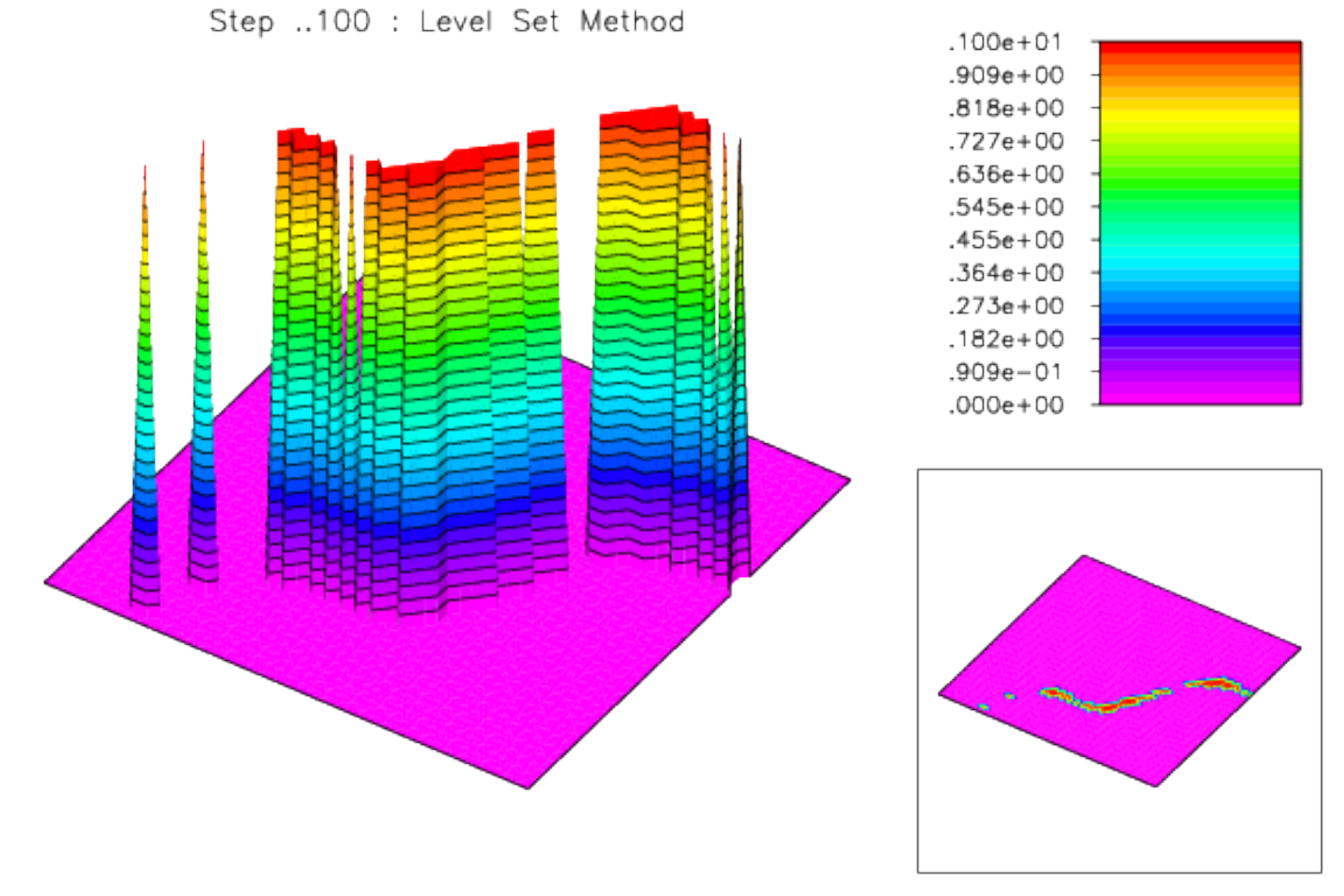} }
\caption{\small Second numerical experiment (analytical P-N junction):
Evolution of the iteration error for the level set method and exact data.}
\label{fig:evol-ls3}
\end{figure}

We conclude this section presenting two different numerical experiments
concerning the identification problem in (\ref{eq:model-DtN}):
\begin{itemize}
\item The first one, for comparison purposes, corresponds to the
identification problem considered in \cite{BELM04} (linear P-N junction;
see Figure~\ref{fig:ls-setup}~(a)). \\
Initially we implemented the level set method for the case of exact data
(see Figure~\ref{fig:evol-ls1}).
Notice that the first picture (top left) corresponds to the initial guess. \\
In a second run we added 10\% random noise to the exact data and repeated
the experiment (see Figure~\ref{fig:evol-ls2}).
\item In the second experiment we try to identify a P-N junction parameterized
by an analytical function (see Figure~\ref{fig:ls-setup}~(b)). Exact data is
used for the reconstruction (see Figure~\ref{fig:evol-ls3}).
\end{itemize}

The voltage-current data pair in the first experiment corresponds to the
the boundary values of the function shown in Figure~\ref{fig:ls-setup}~(c).
This picture shows the solution of (\ref{eq:model-DtN}) for a typical
source $U(x)$.

%-----------------------------------------------------------------------------
\section{Final comments and conclusions} \label{sec:concl}

\subsubsection*{{\em A priori} knowledge of the doping profile}

Due to the choice of $\Omega$, $\partial\Omega_D$ and $\partial\Omega_N$
meet at angles of $\pi/2$. Thus, the solutions of the mixed boundary
value problems in the Landweber-Kaczmarz iteration are not in $H^2(\Omega)$
(see \cite{Gr85} for details). Due to this lack of regularity, the
implementation turns out to be very unstable.
In order to bypass this instability, in \cite{BELM04} the authors made
the additional assumption that the doping profile is known in a thin
strip close to $\partial\Omega_D$. Therefore, only the values of $\gamma(x)$
at a subdomain $\tilde\Omega \subset\subset \Omega$ had to be reconstructed.

Differently from the Landweber-Kaczmarz approach in \cite{BELM04},
the level set method does not require the assumption that the doping profile
is known in some strip close to $\partial\Omega_D$.
For this level set approach, only the knowledge of the doping profile at
$\Gamma_1$ is required in order to obtain a stable performance of the
method.
This weaker assumption agrees with the physical experiment, since we need
to know $\gamma$ at $\Gamma_1$ in order to implement the DtN map in
(\ref{eq:DtN-map}).

\subsubsection*{Amount of data and quality of the reconstruction}

We now comment on the amount of information used in the identification.
In \cite{BELM04} the Landweber-Kaczmarz method was implemented using
different amount of data, i.e., a different number of data voltage-current
pairs. In one of the experiments, a single pair of data was used.
In this case the Landweber-Kaczmarz method reduces to the classical
Landweber iteration.

It is worth noticing that the amount of available data strongly influences
the quality of the reconstruction in the Landweber-Kaczmarz method. However,
observing the results in \cite{BELM04}, no matter how many voltage-current
pairs are available, it does not allow a proper determination of the P-N
junction.

What concerns the quality of the reconstruction of the P-N junction,
the level set approach considered in this paper brings much better
results. In particular if one takes into account that only one pair of
voltage-current data is used.

A possible explanation for the different performance of these methods is
the fact that the Landweber-Kaczmarz approach does not take into account
the assumption that the coefficient $\gamma$ in (\ref{eq:model-DtN}) for
such application is a piecewise constant function.
The Landweber Kazmarz method tries to identify a real function defined on
$\Omega$, which is a much more complicated object than the original unknown
curve (the P-N junction).
Due to the nature of the level set approach, it incorporates in a natural
way the assumption that $\gamma$ is piecewise constant in $\Omega$.

\subsubsection*{Numerical effort}

Next, we compare the numerical effort required for the implementation of
the Land\-weber-Kaczmarz and the level set methods.
If both methods are implemented using a single pair of voltage-current
data, each step of Landweber-Kaczmarz method requires the solution of two
mixed boundary value problems (BVP's), while each level set step requires
the solution of three BVP's.
However, the use of nine pairs of data (as in \cite{BELM04}) requires the
solution of eighteen BVP's in each cycle of the Landweber-Kaczmarz method.
This observation agrees with our numerical tests, where we observed that a
level set step is about ten times faster than a Landweber-Kaczmarz cycle
for nine pairs of voltage-current data.

It is worth noticing that not only the the numerical effort for each step
of the level set method is smaller than the effort for a Landweber-Kaczmarz
cycle, but also the total number of steps required to obtain a good
approximation is smaller than the total number of cycles.
In \cite{BELM04} the iteration was stopped after five thousand cycles (for
both exact and noisy data of 10\%).
This numerical test corresponds to our first experiment, where we needed one
hundred stpes in case of exact data and four hundred steps for noisy data.

\subsection*{Acknowledgments}
The work of A.L. was partially supported by the Brazilian National Research
Council CNPq, grants 305823/03-5 and 478099/04-5.
J.P.Z. acknowledges financial support from CNPq, grants 302161/03-1 and
474085/03-1 and from Prosul program, grant 490300.
P.A.M. acknowledges support from the Austrian National Science Foundation
FWF through his Wittgenstein Award 2000.

%---------------------------------------------------------------------------

\end{document}